\documentclass[acmsmall,natbib=false,nonacm]{acmart}
\AtBeginDocument{%
  }

\RequirePackage[
  datamodel=acmdatamodel,
  style=acmnumeric,%acmauthoryear,
  ]{biblatex}

\usepackage{enumitem}
\usepackage{mathtools}

\definecolor{lightergray}{gray}{0.97}
\definecolor{darkergray}{gray}{0.4}
\usepackage{listings}
\usepackage{menukeys}
\newcommand{\button}[1]{\keys{#1}}
\newcommand{\field}[1]{{\small\textsf{#1}}}
\newcommand{\guielement}[1]{{\small\textsf{#1}}}
\newcommand{\window}[1]{\textsf{#1}}
\newcommand{\odefile}{\textit{odefile}}

\DeclareSourcemap{
  \maps[datatype=bibtex]{
    \map[overwrite]{
      \step[fieldsource=shortjournal, fieldtarget=journaltitle]
      \step[fieldsource=shortseries, fieldtarget=series]
    }
  }
}
\AtEveryBibitem{\clearfield{month}\clearfield{day}}
\newcommand{\D}{\mathop{}\!\mathrm{d}}

\begin{document}

%% BEGIN ACM STANDARD TITLE PAGE

\title{New functionalities in MatCont: delay equations and Lyapunov exponents}


\author{Davide Liessi}
\email{davide.liessi@uniud.it}
\orcid{0000-0002-2649-1253}
\author{Enrico Santi}
\email{santi.enrico@spes.uniud.it}
\orcid{0009-0008-3556-7177}
\author{Rossana Vermiglio}
\email{rossana.vermiglio@uniud.it}
\orcid{0000-0002-7960-1487}
\affiliation{%
  \institution{CDLab -- Computational Dynamics Laboratory, Department of Mathematics, Computer Science and Physics, University of Udine}
  \city{Udine}
  \country{Italy}}

\author{Mayank Thakur}
\email{mamayankt.2002@gmail.com}
\orcid{0009-0004-7973-5719}
\author{Hil G. E. Meijer}
\email{h.g.e.meijer@utwente.nl}
\orcid{0000-0003-1526-3762}
\affiliation{%
  \institution{Department of Applied Mathematics, University of Twente}
  \city{Enschede}
  \country{The Netherlands}}

\author{Francesca Scarabel}
\email{f.scarabel@leeds.ac.uk}
\orcid{0000-0003-0250-4555}
\affiliation{%
  \institution{School of Mathematics, University of Leeds}
  \city{Leeds}
  \country{United Kingdom}}
\affiliation{%
  \institution{CDLab -- Computational Dynamics Laboratory, University of Udine}
  \city{Udine}
  \country{Italy}}

\renewcommand{\shortauthors}{Liessi, Santi, Vermiglio, Thakur, Meijer, Scarabel}

\begin{abstract}
MatCont is a powerful toolbox for numerical bifurcation analysis focussing on smooth ODEs. A user can study equilibria, periodic and connecting orbits, and their stability and bifurcations. Here, we report on additional features in version 7p6. The first is a delay equation importer enabling MatCont users to study a much larger class of models, namely delay equations with finite delay (including delay differential and renewal equations). This importer translates the delay equation into a system of ODEs using a pseudospectral approximation with an order specified by the user. We also implemented Lyapunov exponent computations, event functions for Poincaré maps, and enhanced homoclinic continuation. We demonstrate these features with test cases, such as the Mackey--Glass equation and a renewal equation, and provide additional examples in online tutorials.\end{abstract}

\keywords{Continuation, bifurcation, nonlinear systems, delay differential equations, renewal equations, pseudospectral discretisation, Lyapunov exponents}


\maketitle

%% END ACM STANDARD TITLE PAGE

\section{Introduction}
MatCont is a versatile toolbox for the numerical bifurcation analysis of models emerging from physics, engineering, biology and other fields. This toolbox applies to models formulated as smooth, autonomous ordinary differential equations (ODEs)
\begin{equation}
\label{eq:ode}
x' = f(x,p), \qquad x\in\mathbb{R}^{d}, p\in\mathbb{R}^{m}.
\end{equation}
The user specifies initial data and settings for a set of simulations or continuations. Users can analyse equilibria, limit cycles and connecting orbits, their stability and bifurcations. During and after the computation, the user can inspect and visualise the results.
Since its first release \cite{DhoogeGovaertsKuznetsov2003}, MatCont's functionality has been improved and expanded, but has always focussed on classical bifurcation analysis of ODEs as in \cite{Kuznetsov:2023}, and maps \cite{KuznetsovMeijer:2019}. Surely, there are more classes of models, and methods to characterise their dynamics. For instance, delay equations (DEs) appear naturally in many applications, such as in biology, e.g. population dynamics \cite{Cushing1977} and neuroscience \cite{VisserMeijerVanPuttenVanGils2012}, and physics, e.g. lasers \cite{KrauskopfSieber2023} and climate \cite{FFJin:1997}. Moreover, other dynamical indicators have been increasingly used to characterise the systems' dynamics: for instance, positive Lyapunov exponents (LEs) indicate the appearance of chaos. Here, we present recent changes to MatCont included in \textbf{version 7p6}, enabling users to advance with their model analysis.

Delay equations describe systems in which the current evolution depends not only on the system's current state, but also on its history, via the effect of past events or delayed responses. Here, we focus on autonomous delay differential equations (DDEs) and renewal equations (REs) with finite maximum delay. A typical autonomous DDE with discrete and distributed delays has the form
\begin{equation}\label{dde-proto}
x'(t) = f\left(x(t), x(t-\tau_{1}),\dots,x(t-\tau_{\ell}), \int_{-b}^{-a} g(\theta,x(t+\theta))\D\theta\right),
\end{equation}
while a typical RE is 
\begin{equation}\label{re-proto}
x(t) = f\left(\int_{-b}^{-a} g(\theta,x(t+\theta))\D\theta\right),
\end{equation}
where $\tau_1,\dots,\tau_{\ell}>0$ are discrete delays, $f$ and $g$ may be nonlinear functions depending on parameters, and $0\leq a<b$. We also consider (mixed) systems that couple a finite number of DDEs and REs. In the following, the term DE refers to either type of equation unless otherwise specified.

As the history function is part of the state of a DE, they generate infinite-dimensional dynamical systems. So, progress in the development of numerical tools for DEs has been relatively limited, though several software packages exist, e.g. DDE-BIFTOOL \cite{EngelborghsLuzyaninaRoose2002}, Knut \cite{Szalai2013}, and XPP \cite{Ermentrout2002}. These packages can perform the continuation of branches of equilibria and limit cycles, and detect several bifurcations along solution curves. However, they are restricted to specific classes of models: for instance, DDE-BIFTOOL can treat DDEs with discrete and some state-dependent delays; Knut can treat algebraic equations and neutral DDEs with discrete delays; all of them can only handle discrete delays. Reformulations of the systems allow treating distributed delays and REs in an ad hoc manner. However, the learning curve of these toolboxes is steep, while a user-friendly interface is missing.

One can draw upon the well-established software tools for ODEs by several methods, such as the linear chain trick for discrete delays \cite[Chapter 7]{Smith2011} (see also \cite{AndoBredaGava2020}) or Volterra kernels for distributed delays \cite{NevermannGros2023}. These methods, however, are either not general or provide poor approximations. Another more recent pseudospectral method \cite{BredaDiekmannGyllenbergScarabelVermiglio2016,ScarabelDiekmannVermiglio2021} uses the formulation of the DE as a semilinear abstract differential equation to obtain an approximating ODE via collocation on optimal nodes. This approach suggests how to treat the time shift of the history function and its discretisation, and has several advantages. First, it can treat general DEs. Second, the approximation provides good results with low orders, and the method has a variable order and allows to assess the accuracy of the numerical results. Here, we describe how this ODE approximation is implemented in MatCont using a \textit{delay equation importer}. This importer expands the class of models that can be studied with a single software package, saving time and learning effort for modellers.

Another addition to MatCont is the seamless integration for computing LEs. This method has been around for a long time already, see, e.g., \cite{Govorukhin:2004}. Remarkably, modellers always have had to do this analysis ad hoc, while MatCont offers precisely the necessary elements, i.e. the system plus the Jacobian matrix for the variational equation.  

Next to these additions, we created scripts for automated testing of tutorials for new versions and event functions for Poincar\'e maps.

In this paper, we describe our design choices for the DE importer so that the resulting functionality builds on rather than mixes with the ODE approach of MatCont. Indeed, the MatCont user works with ODEs, while several aspects of the implementation are DE-specific and are not relevant to the usual analysis of ODEs. Also, for the computation of LEs our design follows the idea of first specifying settings before the actual computation.
We also provide guidelines for these numerical settings. We report on performance in run-time and approximation error for two examples, while our online documentation has more tutorials.

\section{The Delay Equation Importer}

Sun-star calculus provides an abstract framework for DEs to establish the fundamental results of dynamical system theory, such as the principle of linearised stability, theorems on local invariant manifolds and bifurcations, and Floquet theory for periodic solutions \cite{BredaLiessi2021,DiekmannGettoGyllenberg2008,DiekmannVanGilsVerduynLunelWalther1995}. In this framework, a DE is reformulated as a semilinear abstract differential equation. In this way, the evolution law of the dynamical system describes both the change of the state and the shift of the history.\footnote{This is known as shift-and-extend (into the future).} This idea requires a larger state space as the history function has to be included as well. Such history function is numerically represented with collocation polynomials of degree $M$ defining the order of the approximation \cite{BredaDiekmannGyllenbergScarabelVermiglio2016}. Chebyshev points of the second kind were chosen as to cover the endpoint at the maximal delay. The action of the shift is differentiation, which is linear, minimising additional complexity in formulating the evolution on the extended approximating state space. This setup allows for formulating an ODE that approximates the DE, with $M$ auxiliary variables for each model variable.
These choices give a pseudospectral method that approximates the eigenvalues associated with equilibria with exponential order of convergence in $M$ \cite{BredaMasetVermiglio2015,DiekmannScarabelVermiglio2020,ScarabelDiekmannVermiglio2021,DeWolffScarabelVerduynLunelDiekmann2021}. We note that the convergence rate is an open problem for limit cycles. Tests of this approximation method have already been performed with MatCont since \cite{BredaDiekmannGyllenbergScarabelVermiglio2016}.
Appendix \ref{sec:discretization} summarises the main aspects of the reformulation and discretisation approach.

The implementation of the approximating ODE and its analysis via command line have
so far been a manual task, which constituted a barrier to a wider adoption of this method. In MatCont, a user can enter an ODE system using the graphical interface (GUI), and the input is converted to a MATLAB file. We created a \textit{delay equation importer} so that a user can declare a model of a new class, allowing additional syntax to specify time-delayed arguments and integrals. Once the user has specified the DE, the syntax is interpreted, and code for the pseudospectral approximation is generated, adding the appropriate number of auxiliary variables. To improve the user experience while studying the approximating ODE within the GUI, we have made several changes in how MatCont handles particular data, which are described below.
All the changes have been guided by two principles: on the one hand, we strived to be transparent to the user, maintaining the auxiliary variables relevant to the DE both accessible and visible; on the other hand, we facilitated where possible the reconstruction and plot of the original model states via the GUI.

\paragraph{Importer}

\begin{figure}
    \centering
    \includegraphics[width=0.6\linewidth]{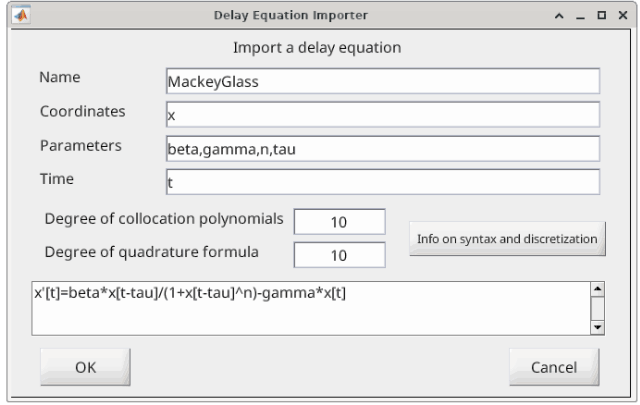}
    \Description{Screenshot of MatCont's \window{Delay Equation Importer} window (see the caption).}
    \caption{The \window{Delay Equation Importer} window filled out with the input for the Mackey--Glass DDE \eqref{Mackey-Glass}.}
    \label{fig:importer}
\end{figure}

The \window{Delay Equation Importer} window can be opened from MatCont's main window by selecting \menu{Select>System>Delay Equation Importer>New} (see Figure~\ref{fig:importer}).
The importer window translates a DE into the appropriate system definition file (or \odefile{} in short) to be used for the analyses in MatCont.
The \odefile{} contains the definition of the ODE obtained via pseudospectral approximation (as described in Appendix \ref{sec:discretization}) by automatically defining all the required discretisation parameters and matrices. 
Systems can be modified via \menu{Select>System>Delay Equation Importer>Edit Current System}; doing it
with the standard editor for ODEs will usually result in an error.

The degree $M$ of the collocation polynomials and the degree of the quadrature formula have to be specified in the importer, since different values give different approximating systems with incompatible representations of the results.
Hence, to experiment with different values of the numerical parameters, the user should define a new system.
However, one can take advantage of the previous definition via \menu{Select>System>Delay Equation Importer>Edit Current System}. Changing also the \field{Name} field will result in the creation of a new \odefile{} under the new name.

\paragraph{Input syntax}
MatCont's syntax for ODEs is based on MATLAB's syntax and aims at mirroring the usual mathematical notation; e.g., the ODE $x' = x(1-x) - \frac{2xy}{a+x}$ is input as \lstinline!x'=x*(1-x)-2*x*y/(a+x)!.
For the importer, we extend that syntax to DEs with the following elements.

Time dependency of the unknown functions is specified with square brackets: e.g., $x(t-1)$ is input as \lstinline!x[t-1]!.
The argument of the square bracket must be in the form \lstinline!t-!\emph{expr} (or \lstinline!t+!\emph{expr}), with the time variable \lstinline!t! as first character, where \lstinline!-!\emph{expr} (or, respectively, \lstinline!+!\emph{expr}) is an expression resulting in a negative number.
Specifying the current time is also allowed, so the expressions \lstinline!x! and \lstinline!x[t]! are equivalent.
Integrals for distributed delays use a function-like
syntax: $\int_{-b}^{-a} g(\theta,x(t+\theta))\D\theta$ is input as \lstinline!DE_int(@(theta)g(theta,x[t+theta]),-b,-a)!.%
\footnote{In fact, \lstinline!DE_int! is a local function defined in the \odefile.}
The discrete delays and the endpoints of distributed delays (\lstinline!a! and \lstinline!b! above) may be numbers or expressions, possibly depending on parameters.

Like in the case of ODEs, the derivative of unknown functions on the left-hand side of DDEs is marked by a prime sign, as \lstinline!x'! or, equivalently, \lstinline!x'[t]!; for REs, instead, the left-hand side consists of the unknown function itself, specified with \lstinline!x! or \lstinline!x[t]!.

As an example, equations \eqref{dde-proto} and \eqref{re-proto} can be written, respectively, as
\begin{lstlisting}
x'[t]=f(x[t],x[t-tau_1],...,x[t-tau_l],DE_int(@(theta)g(theta,x[t+theta]),-b,-a))
\end{lstlisting}
and
\begin{lstlisting}
x[t]=f(DE_int(@(theta)g(theta,x[t+theta]),-b,-a))
\end{lstlisting}
where \lstinline!tau_1!, \dots, \lstinline!tau_l!, \lstinline!a!, and \lstinline!b! should be defined as parameters of the system, and \lstinline!f! and \lstinline!g! should be valid MATLAB functions.
A guide to the input syntax can be accessed from the \window{Delay Equation Importer} window via the button \button{Info on syntax and discretization}.

Section~\ref{sec:examples} contains more examples and the corresponding syntax.
In particular, examples~\eqref{eq:logdaphnia} and~\eqref{eq:neuron2} show that, as in MatCont's ODE input, the user can define intermediate variables and use them in the right-hand side of the equations.

\paragraph{Sanity checks}
The \odefile{} includes a check on the delays of the system, which is essential since delays can depend on the parameters.
The continuation is stopped whenever a delay becomes negative or the maximal delay is zero.

\paragraph{Auxiliary variables, labelling, and reconstruction of the solution of renewal equations}
For each variable of the DE, the approximation method introduces auxiliary variables which describe the history function (for  DDEs) or a primitive of the history function (for REs); see Appendix~\ref{sec:discretization} for details.
The $k$-th auxiliary variable related to the $i$-th unknown is labelled (e.g. in the \window{Starter} or \window{Numeric} window) with the name of the original unknown and the suffix \lstinline!_aux!$k$ (e.g., \lstinline!x_aux10!).\footnote{Using the notation of Appendix \ref{sec:discretization}, this variable corresponds to the $i$-th component of the vector $V_{M,k}$.}

\begin{figure}
    \centering
    \raisebox{-0.5\height}{\includegraphics[width=0.6\linewidth]{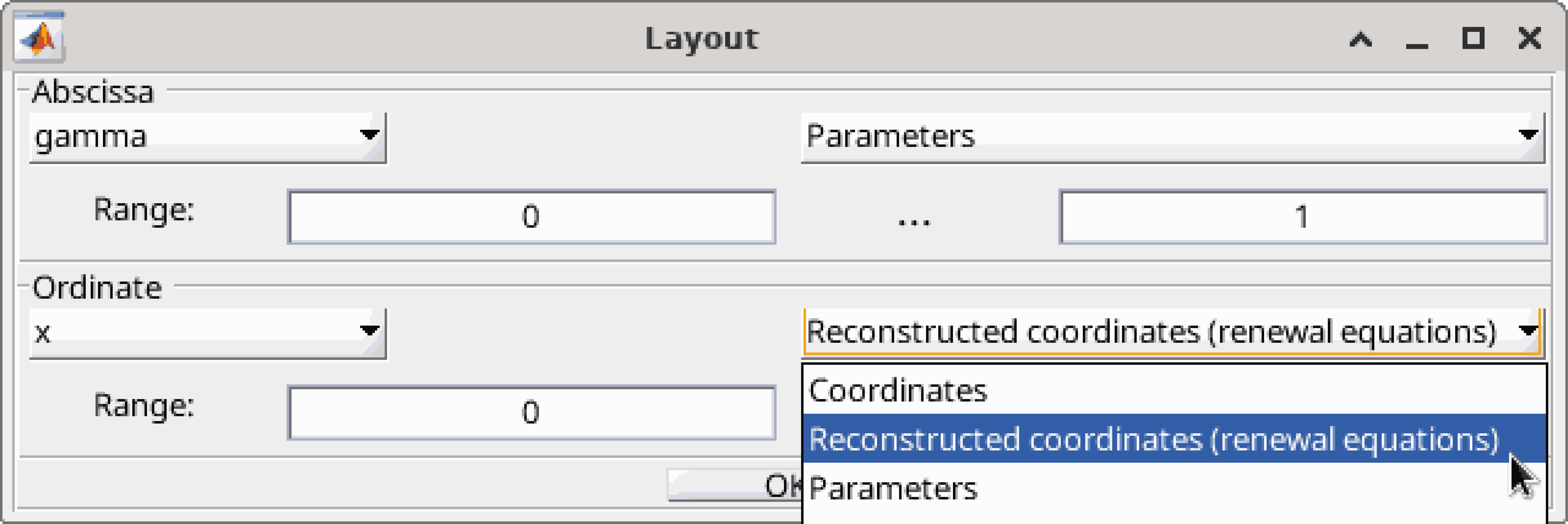}}%
    \hfill%
    \raisebox{-0.5\height}{\includegraphics[width=0.37\linewidth]{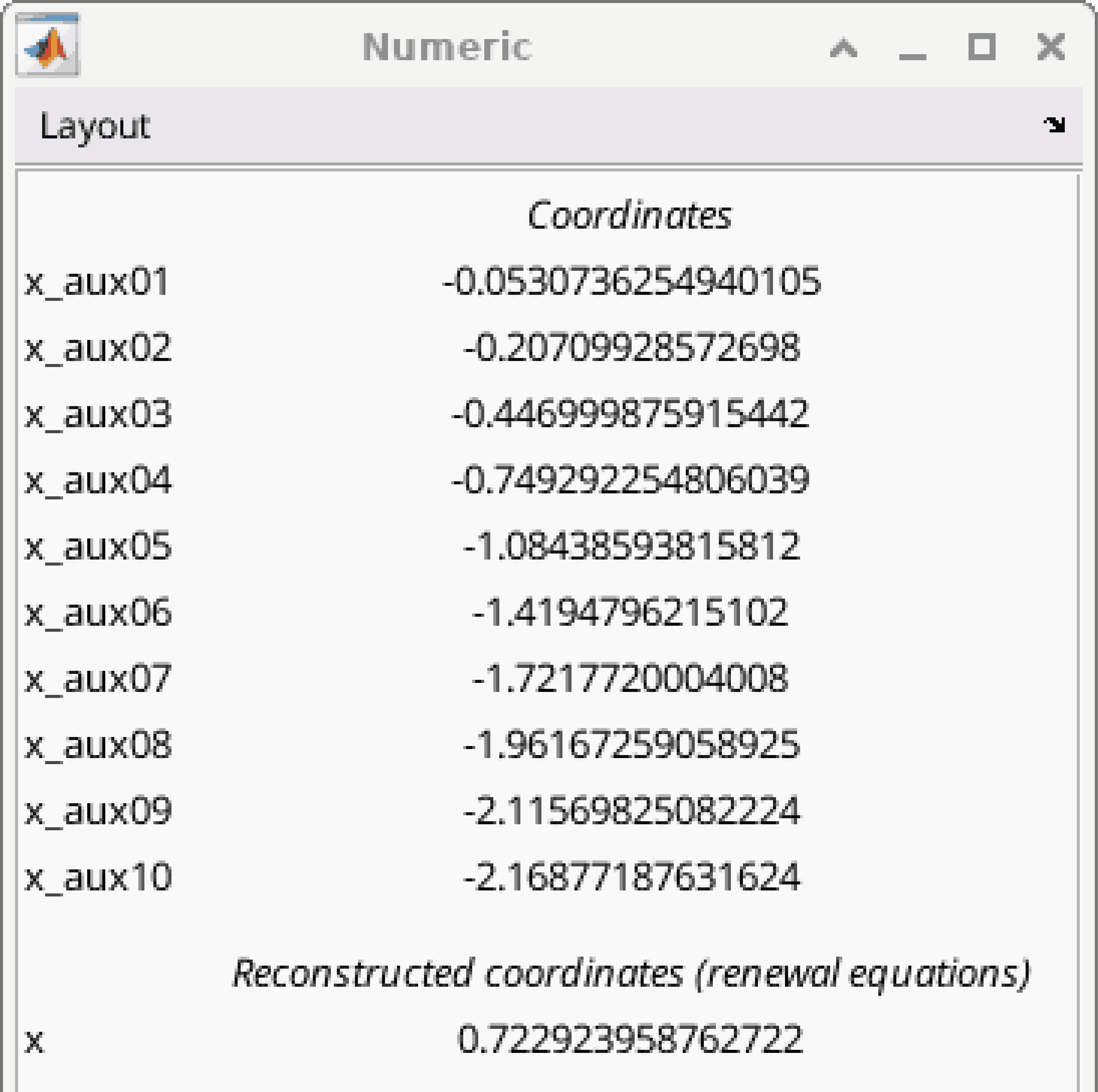}}%
    \Description{Screenshots of MatCont's \window{Layout} and \window{Numeric} windows (see the caption).}
    \caption{Layout settings for a graphical window, showing the reconstructed coordinates for an RE (left), and the \window{Numeric} window for an RE, showing both the auxiliary and the reconstructed coordinates.}
    \label{fig:re-plot-properties-and-numeric}
\end{figure}

In the case of DDEs, the first components of the approximating ODE system approximate the values of each unknown function at the current time $t$, and are thus labelled with the name of the original variables (e.g., \lstinline!x!).\footnote{Using the notation of Appendix \ref{sec:discretization}, these variables correspond to the entries of $x_M$.}
In the case of REs, instead, all variables of the approximating ODE system are auxiliary,
and the solution must be reconstructed using~\eqref{eq:re-extra-step}.
The reconstructed solution is made available for plotting and monitoring in the appropriate windows (see Figure~\ref{fig:re-plot-properties-and-numeric}), thanks to an adaptation to the corresponding files in MatCont's main codebase.

\begin{figure}
\centering
\includegraphics[clip,trim=20mm 157mm 17mm 30mm]{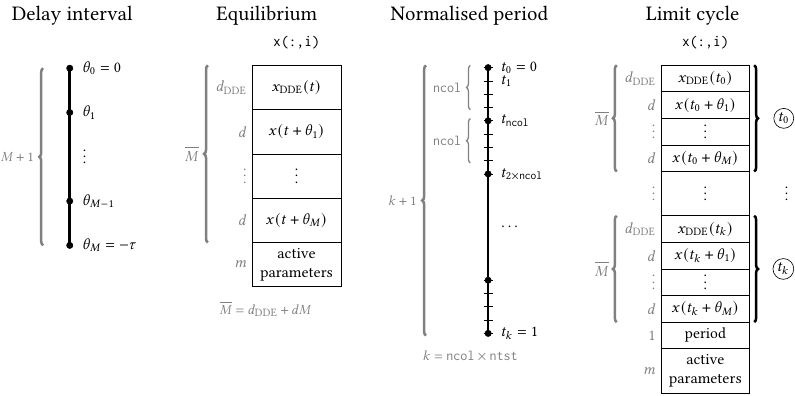}
\Description{Schemes of the state vectors of the system of ODEs approximating a system of DEs (see the caption).}
\caption{State vectors of the system of ODEs approximating a system of DEs (DDEs, REs, or coupled DDEs/REs) with $d$ coordinates.
We denote by $d_{\text{DDE}}$ the number of coordinates described by DDEs, while the remaining $d-d_{\text{DDE}}$ coordinates are described by REs.
For a system of DDEs only, $d=d_{\text{DDE}}$; for a system of REs only, $d_{\text{DDE}}=0$.
The diagrams show the notation for the collocation of the history functions in the delay interval (first from the left) and of the limit cycles in one period interval normalised to $[0,1]$ (third), and the structure of the columns (one for each continuation point) of the \lstinline!x! output variable of the \lstinline!cont! function, for equilibria (second) and limit cycles (fourth).
}
    \label{fig:limit-cycle}
\end{figure}

Figure~\ref{fig:limit-cycle} shows schematic representations of the state vector of the approximating ODE, which can be a useful reference when the user needs to access the data for ad-hoc analyses or visualisations.

\paragraph{Setting the initial point}
Both the time integration of the equation and the continuation of an equilibrium branch require setting an initial point.
For DEs, this corresponds to an initial function on the delay interval, which has to be computed at the collocation nodes and, in the case of REs, integrated via \eqref{eq:re-initial-value} to construct the initial value for the approximating ODE.

\begin{figure}
    \centering
    \includegraphics[height=0.3\linewidth]{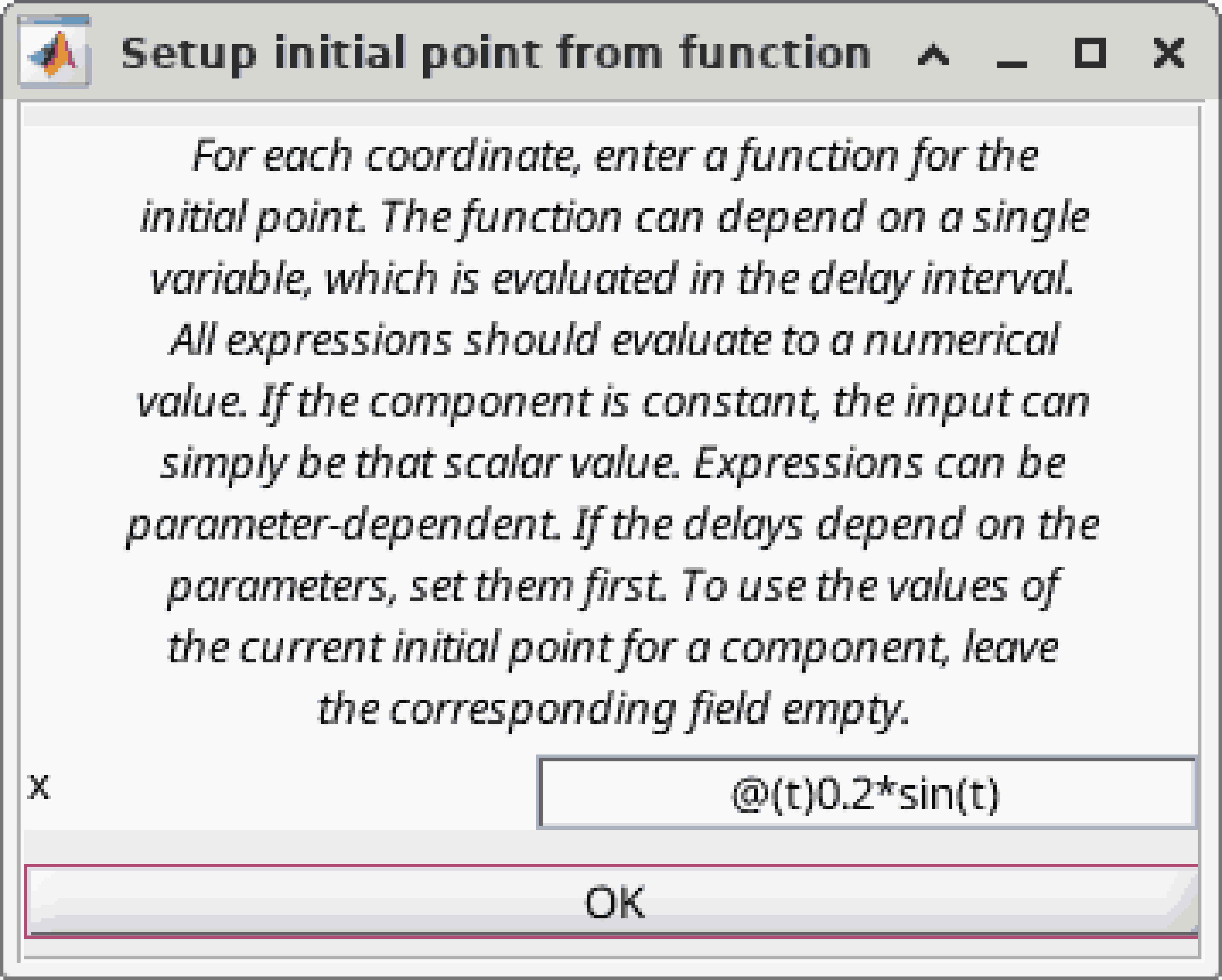}%
    \qquad%
    \includegraphics[height=0.3\linewidth]{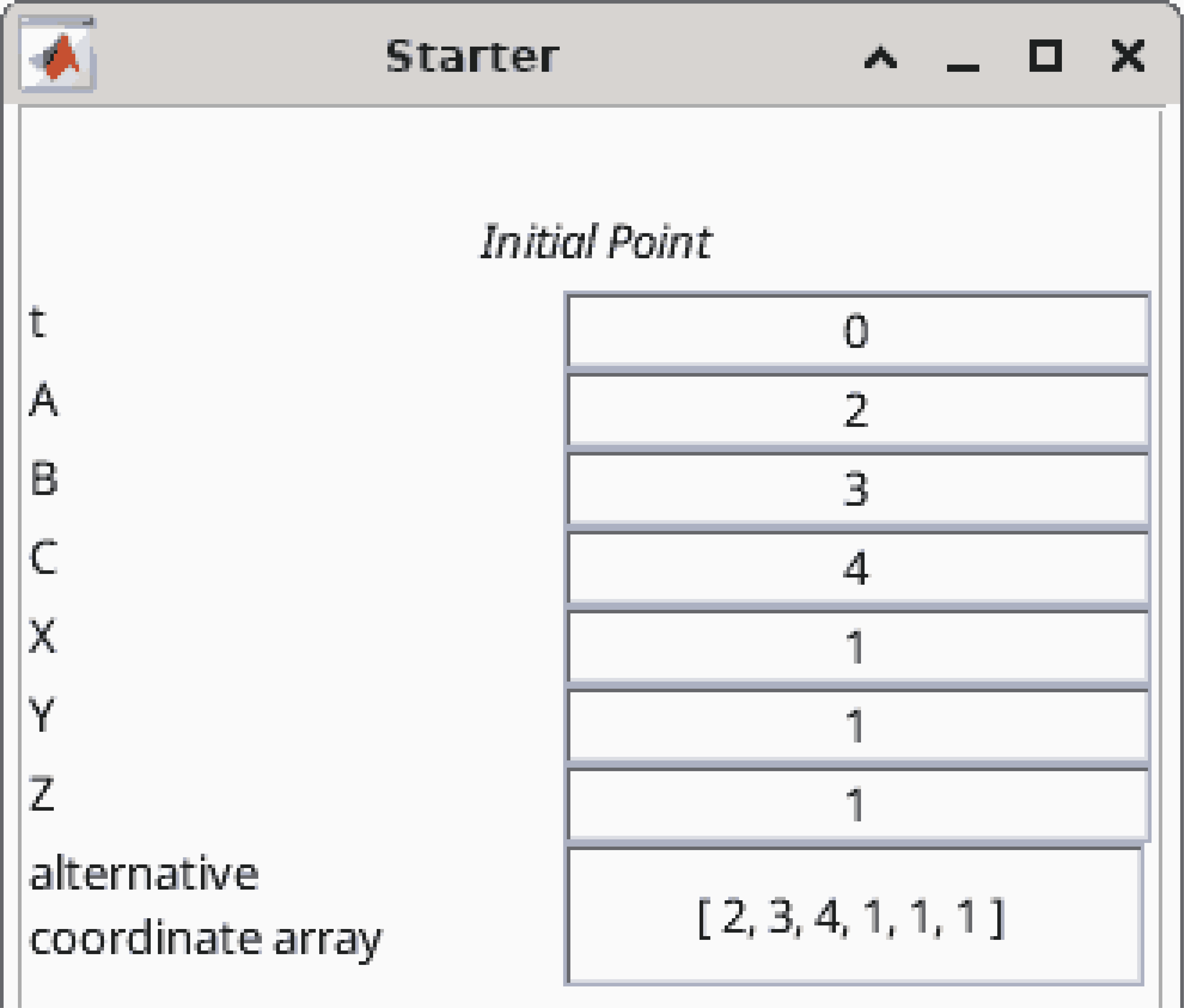}
    \Description{Screenshots of MatCont's \window{Setup initial point from function} and \window{Starter} windows (see the caption).}
    \caption{The window to setup the initial point from a function (left), and the \window{Starter} window with the new \field{coordinate array} field, showing the result of inserting \texttt{[2:4,ones(1,3)]} in that field for a system of six ODEs in the variables $A$, $B$, $C$, $X$, $Y$, and $Z$ (right).}
    \label{fig:initial-value-tool-and-vector-coordinates}
\end{figure}

To do this, we added the button \button{Setup initial point from function} to the \window{Starter} window.
The button opens a window (Figure~\ref{fig:initial-value-tool-and-vector-coordinates}, left) where the user can specify a function handle (or a number if the initial function is constant) for each variable of the DE.
The setup window then calls a function that correctly computes and loads the initial vector for the ODE.
This function is provided in the \odefile,
and can thus also be accessed when using MatCont from the command line.

\paragraph{Numerical choices}

As collocation nodes in the maximal delay interval $[-\tau,0]$, we chose the Chebyshev points of the second kind.
Interpolation is done with the barycentric Lagrange interpolation formula \cite{BerrutTrefethen2004}.
For the quadrature of integrals we chose the Clenshaw--Curtis formula \cite{ClenshawCurtis1960,Trefethen2008}.
For more details, see Appendix~\ref{sec:discretization}.

\paragraph{Computation of the Jacobian matrix}
Thanks to the semilinear structure of the approximating ODE (see \eqref{eq:dde-discr} and \eqref{eq:re-discr}), a large part of the Jacobian matrix can be explicitly defined independently of the right-hand side of the DE.
To exploit this property, we define the first-order derivative using numerical differentiation for the nonlinear part only.
The computation of the derivatives via user-specified functions or via the Symbolic Math Toolbox is disabled for imported DEs.

\section{Lyapunov exponents and other changes}

\paragraph{Lyapunov exponents}
The implementation employs the classical QR algorithm, as provided by Govorukhin on MATLAB Central File Exchange (with permission) \cite{Govorukhin:2004}. For this QR method, one integrates the variational equation $Y'=Df(x(t),p)Y$ with initial condition $Y(0)=I_{d}$ along an orbit $x(t)$ for \eqref{eq:ode}. After fixed intervals, one collects the growth factors and re-orthogonalises the Lyapunov vectors represented by the columns of $Y$. These growth factors are used to compute all LEs. Once a user has entered a system in MatCont, it is straightforward to compute the Jacobian matrix $Df$ along an orbit, and to integrate the variational equation. We have added a curve type \guielement{LyapunovExponent (LyapExp)} that calls this numerical method, avoiding the need to create separate codes.

The user can choose the curve type after selecting the initial point. For LEs, the initial point has no particular properties, so the user can choose any point, including the type \guielement{Point} labelled by \texttt{P}. Next, the user specifies the type of a solution curve: instead of \guielement{Orbit (O)}, typical for an arbitrary point, one chooses the new curve type \guielement{LyapunovExponent (LyapExp)}. We provide two options. First, for fixed values of the parameters, MatCont computes the running estimate of LEs. With this option, the user can check numerical settings and inspect the convergence of LEs. Second, the user can specify an array of values for a single component of the parameter vector $p$. For each of these values, the LEs of an attractor are computed, and the last state can be used for the next value. This option allows to see how LEs evolve as the parameter changes and also see bifurcations and indicate the transition to chaos. The user can also call the method from the command line, which gives more control.

\paragraph{Setting the initial point}
To simplify the definition of the initial point in the \window{Starter} window in case of ODEs of large dimension, we added a field called \field{coordinate array} (Figure~\ref{fig:initial-value-tool-and-vector-coordinates}, right) which allows the user to specify the values of the coordinates as an array (in any valid MATLAB notation, including via function calls).
This new field and the fields for the individual coordinates are automatically kept synchronised.

\paragraph{Other improvements}
Finally, we mention that we have integrated a new predictor for the homoclinic bifurcation curve presented in \cite{Bosschaert:2024}. The algorithm is the culmination of a series of papers on higher-order approximations of this orbit on the parameter-dependent centre manifold. Motivated by recent work on high-dimensional ODE models \cite{Pusuluri:2021,Scully:2025}, we improved the test functions for homoclinic bifurcations.

Eigenvalues and multipliers are now sorted by decreasing relevance (i.e., real part for eigenvalues and magnitude for multipliers). A new option in the \window{Numeric} window allows to show only a certain (configurable) number of the most relevant eigenvalues.
This option is available via the \menu{Options>Eigenvalues/Multipliers} menu.
These changes are particularly useful for imported DEs, and possibly discretised PDEs, as the dimension of the resulting ODE system can be large.

\section{Illustrative examples}
\label{sec:examples}

In the following we collect some examples to illustrate the new capabilities of MatCont.
We present a numerical bifurcation analysis of a DDE and we demonstrate the computation of LEs of an ODE and an RE.
We also collect some more examples of DEs to show different elements of the syntax.

\subsection{The Mackey--Glass delay differential equation}

We here present a numerical bifurcation analysis of the well-known Mackey--Glass DDE \cite{MackeyGlass1977,BredaDiekmannGyllenbergScarabelVermiglio2016}, including step-by-step instructions to reproduce the analysis in MatCont via the GUI. 
We refer to the PDF and video tutorials available at \url{https://sourceforge.net/projects/matcont/} for a more complete and extensive demonstration.

The model reads as
\begin{equation}\label{Mackey-Glass}
x'(t) = \beta \frac{x(t-\tau)}{1+x(t-\tau)^n} - \gamma x(t),
\end{equation}
for positive parameters $\beta$, $\gamma$ and $n$, and one single discrete delay $\tau>0$.

For this demonstration, we open the \window{Delay Equation Importer} window from MatCont's main window (titled \window{MatCont GUI}) via 
\menu{Select>System>Delay Equation Importer>New}, and we input the model parameters and the equation (see Figure \ref{fig:importer}). 
The input syntax for equation \eqref{Mackey-Glass} is 
\begin{lstlisting}
x'[t]=beta*x[t-tau]/(1+x[t-tau]^n)-gamma*x[t]
\end{lstlisting}
with \lstinline!x! as the list of coordinates and \lstinline!beta,gamma,n,tau! as the list of parameters.
In this test, we use the default degree $10$ for the collocation polynomials (the quadrature formula is not relevant in this example). 
Pressing \button{OK} creates the \odefile{} in the \directory{Systems} folder (\directory{MackeyGlass.m} in our case) and loads the system.
The \odefile{} contains details about the different sections of the code.

\begin{figure}
    \centering
    \includegraphics[width=0.4\linewidth]{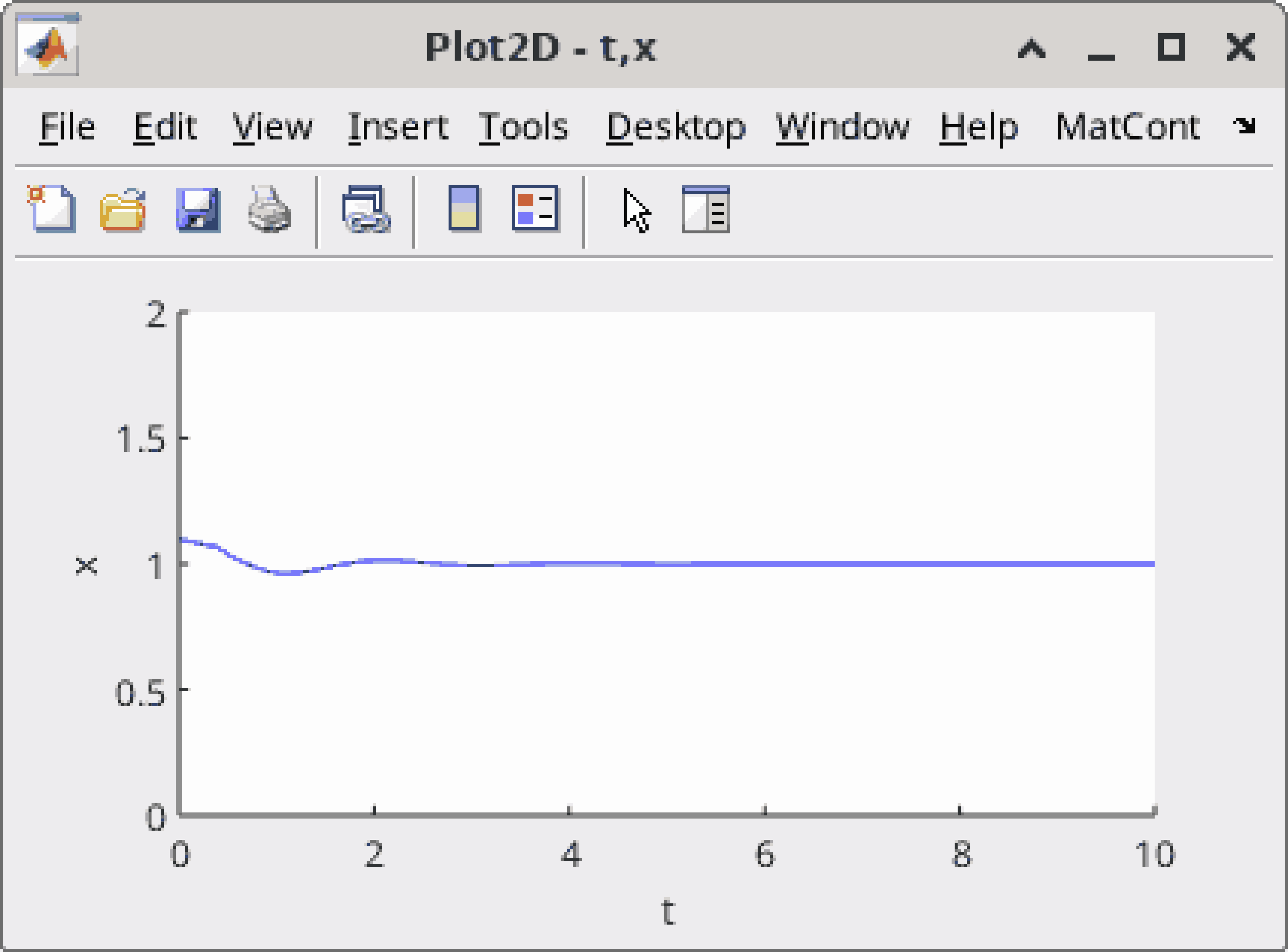}%
    \qquad%
    \includegraphics[width=0.4\linewidth]{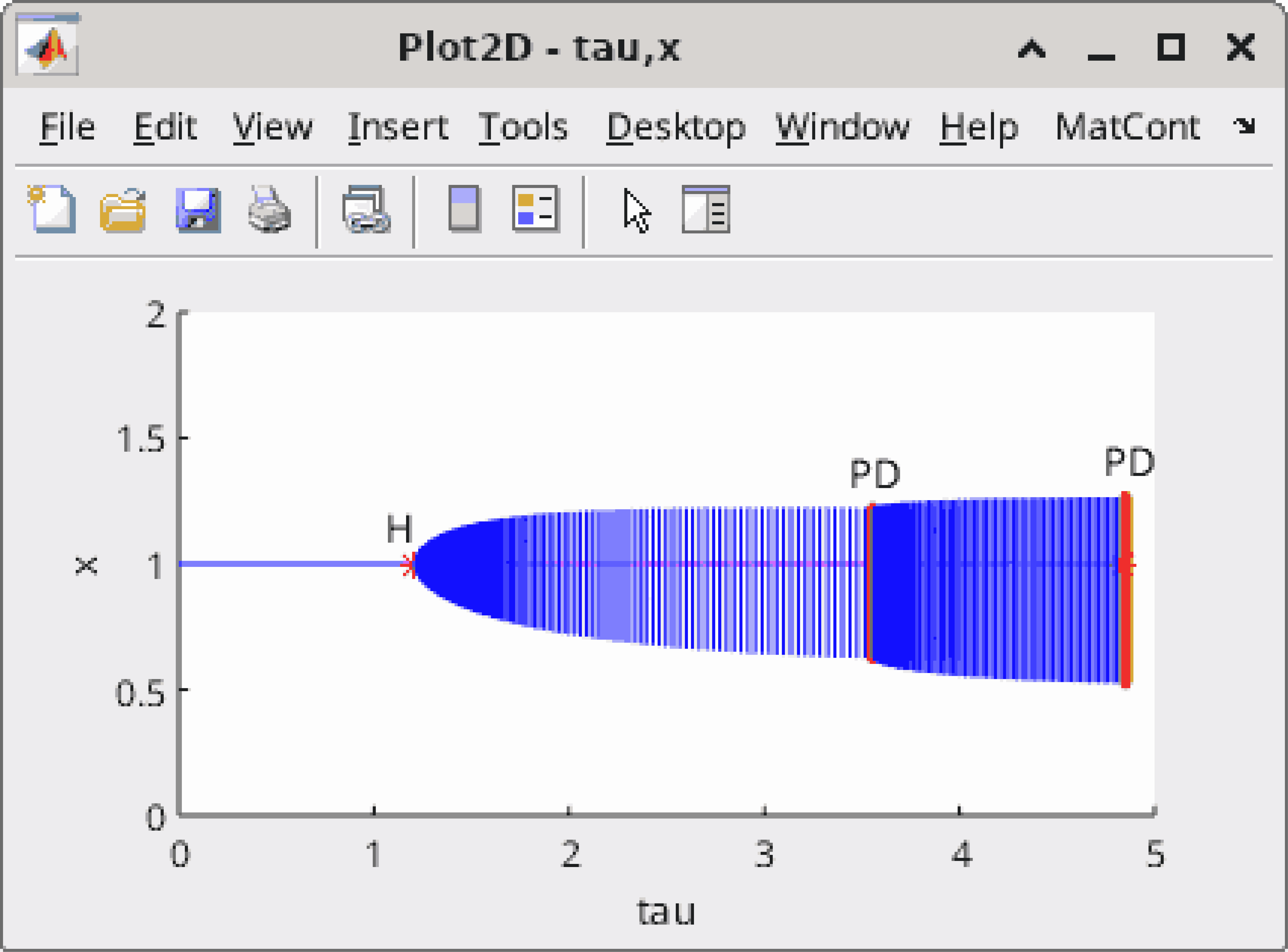}
    \Description{Screenshots of MatCont's \window{Plot2D} windows containing the plots described in the caption.}
    \caption{Mackey--Glass equation: the orbit obtained from the initial vector \lstinline![1.1,1.05,1.2,ones(1,8)]! (left) and the bifurcation diagram with respect to $\tau$ (right).}
    \label{fig:mg-triv-eq-and-bif}
\end{figure}

We start a time integration from a perturbation of the equilibrium $x=1$.
To do this, from the main window we select \menu{Type>Initial Point>Point}.
In the \window{Starter} window we set the initial value to the one in the caption of Figure \ref{fig:mg-triv-eq-and-bif}, and the parameters $\beta=2$, $\gamma=1$, $n=6$ and $\tau=0.5$.
In the \window{Integrator} window we set \field{Interval} to $10$.
We open a graphical window with \menu{Window/Output>Graphic>2D plot} and select the time variable $t$ as abscissa with range $[0,10]$ and the coordinate $x$ as ordinate with range $[0,2]$ via \menu{MatCont>Layout}.
From the main window we select \menu{Compute>Forward}. 
As expected, the orbit seems to converge to an equilibrium (see the left panel in Figure \ref{fig:mg-triv-eq-and-bif}).

To perform a continuation of that equilibrium, we select the last point of the time integration as starting point: we open the \window{Data Browser} window by clicking on \button{View Result} in the \window{Control} window, then select \guielement{Last Point} and click on \button{Select Point}.
This operation loads the values in the \window{Starter} window.
We clear the graphical window with \menu{MatCont>Clear} and set the parameter $\tau$ as abscissa with range $[0,5]$.
To plot stable and unstable branches in different colors, we select \menu{MatCont>Plot Properties} and set the style \lstinline!'Color','magenta','linestyle','-'! in the fields \field{EP if unstable} and \field{LC if unstable}.
In the main window we select \menu{Type>Initial Point>Equilibrium}, then \menu{Compute>Backward}.
The continuation proceeds in the direction of decreasing $\tau$ and stops with the message \lstinline!Current stepsize too small!.
This happens because the delay $\tau$ is approaching $0$: the continuation step lands on a point with negative delay, where the computation cannot continue, so MatCont tries to reduce the step size until it becomes smaller than \field{MinStepSize} (see the \window{Continuer} window).
Indeed, the value of $\tau$ at the last point of the continuation curve is approximately $1.2\text{e-}5$, as we can see in the \window{Data Browser}.
We now select \menu{Compute>Forward} and the continuation detects a Hopf point at $\tau \approx 1.21$. We select \button{Resume} and then stop the continuation after a few more points, in order to observe the change of stability in the plot. In the \window{Data Browser} we can see that the first Lyapunov coefficient corresponding to the first Hopf bifurcation is negative, suggesting a supercritical bifurcation. 

From the Hopf bifurcation we initialise the computation of a branch of limit cycles. 
In the \window{Data Browser} we select the Hopf point and load it in the \window{Starter} window. Note that, in the main window, the \field{Initial Point Type} has changed to \guielement{Hopf (H)} and the \field{Curve Type} to \guielement{Limit cycle (LC)}. 
Observe also that in the \window{Starter} window a new parameter \field{Period} has appeared, which is automatically selected for the continuation of the limit cycles. 
The window also includes the standard numerical parameters \field{ntst} and \field{ncol} for the approximation of the limit cycles. 
We start the continuation by \menu{Compute>Forward} and resume it if necessary until we find a period-doubling bifurcation at $\tau \approx 3.56$, extending the continuation for a few more points. We observe that the limit cycles after the bifurcation are unstable (we can verify in the \window{Data Browser} that \field{pdcoefficient} is negative).

We now select the \guielement{PD} point as initial point to compute the branch of double-period orbits emerging here. 
The continuation 
detects another period doubling at $\tau \approx 4.86$.
Figure~\ref{fig:mg-triv-eq-and-bif} (right panel) shows the bifurcation diagram computed so far.

\subsubsection{Convergence and comparison with other methods}
\label{sec:mgconvergence}

The accuracy of the continuation and bifurcation detection depends on the order of the approximation, and in particular a larger order may be required to detect more complex bifurcations. 

Table~\ref{tab:mg} shows a comparison of the output obtained with different orders (`pseudospectral', $M=5,10,20,40$) and methods (`linear chain trick', `DDE-BIFTOOL'). 
We present the computed values of the two Hopf points on the nontrivial equilibrium branch and the two period doublings, and the computation times to perform the continuation of the branches of equilibria, limit cycles, and doubled cycles. 
The values are obtained with accuracy options $10^{-8}$, for 100 continuation points.

In this example, the pseudospectral method implemented within MatCont is able to capture the bifurcations reasonably well, while the computation times are generally higher than those of DDE-BIFTOOL, which is a package developed specifically for DDE with discrete delays. 
We note that an exact quantitative comparison is hard to achieve because of the different methods and design choices at the core of the different software packages. 
The codes used to produce the table are available in MatCont's \directory{Testruns/TestDelayMG} folder.

\begin{table}
    \centering
    \begin{tabular}{rcccc}
    \toprule
    $M$ & 1\textsuperscript{st} Hopf & 2\textsuperscript{nd} Hopf & 1\textsuperscript{st} PD & 2\textsuperscript{nd} PD \\
    \midrule
    linear chain trick: 10 & 1.98041101 & --- & --- & --- \\
    20 & 1.48084288 & --- & --- & --- \\
    40 & 1.32851441 & --- & 5.75976022 & --- \\
    80 & 1.26550447 & --- & 4.27685322 & 6.85489462 \\
    \midrule
    pseudospectral: \ 5 & 1.20947760 & --- & 3.48977702 & 4.88889653 \\
    10 & 1.20919957 & 4.83374221 & 3.55626313 & 4.85243716 \\
    20 & 1.20919957 & 4.83679830 & 3.55716811 & 4.86050005 \\
    40 & 1.20919957 & 4.83679830 & 3.55716813 & 4.86049811 \\
    \midrule    
    DDE-BIFTOOL & 1.20919957 & 4.83679830 & 3.55716847 & 4.86050131 \\
    \bottomrule \\
    \toprule
    \textit{Computation times ($s$):} & & Equilibrium & Limit cycle & Doubled cycle \\
    \midrule
    linear chain trick: 10 & & 0.2 & 21.1 & --- \\
    20 & & 0.7 & 52.7 & --- \\
    40 & & 2.7 & 1030.8 & 242.3 \\
    80 & & 108.0 & 12455.3 & 12453.2 \\
    \midrule
    pseudospectral: \ 5 & & 0.7 & 11.8 & 11.9 \\
    10 & & 0.2 & 31.0 & 31.0 \\
    20 & & 0.7 & 114.5 & 94.4 \\
    40 & & 4.7 & 484.7 & 371.3 \\
    \midrule
    DDE-BIFTOOL & & 7.0 & 11.2 & 13.7 \\
    \bottomrule
    \end{tabular}
    \smallskip
    \caption{Mackey--Glass equation: approximated values of $\tau$ at the bifurcation points for different values of $M$ (pseudospectral) and different methods (linear chain trick, DDE-BIFTOOL).
    The second Hopf bifurcation appears in the unstable part of the nontrivial equilibrium branch. The lower part of the table contains the computation time to perform the continuation of the branch of equilibria, limit cycles and doubled cycles. We set the maximum number of points in the branch to $100$, and we set the tolerances of MatCont and DDE-BIFTOOL at $10^{-8}$. The experiments were done with MATLAB 2023a on an Apple MacBook Air (M1, 2020) with 8\,GB of RAM running macOS Sonoma 14.7.2.
    }
    \label{tab:mg}
\end{table}

\subsection{Examples for Lyapunov exponents}
\begin{figure}
    \centering
    \includegraphics[width=0.32\linewidth]{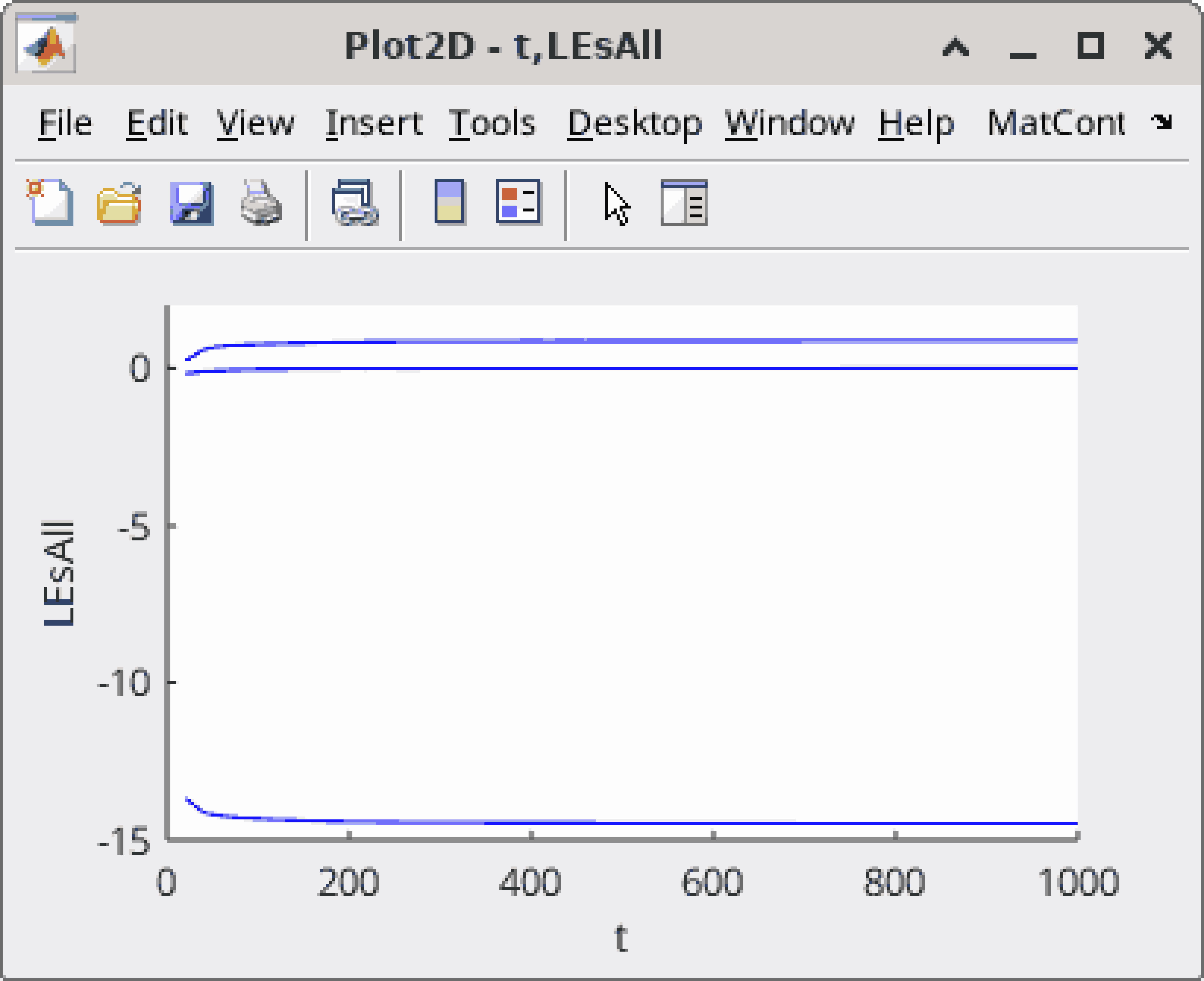}%
    \hfill%
    \includegraphics[width=0.32\linewidth]{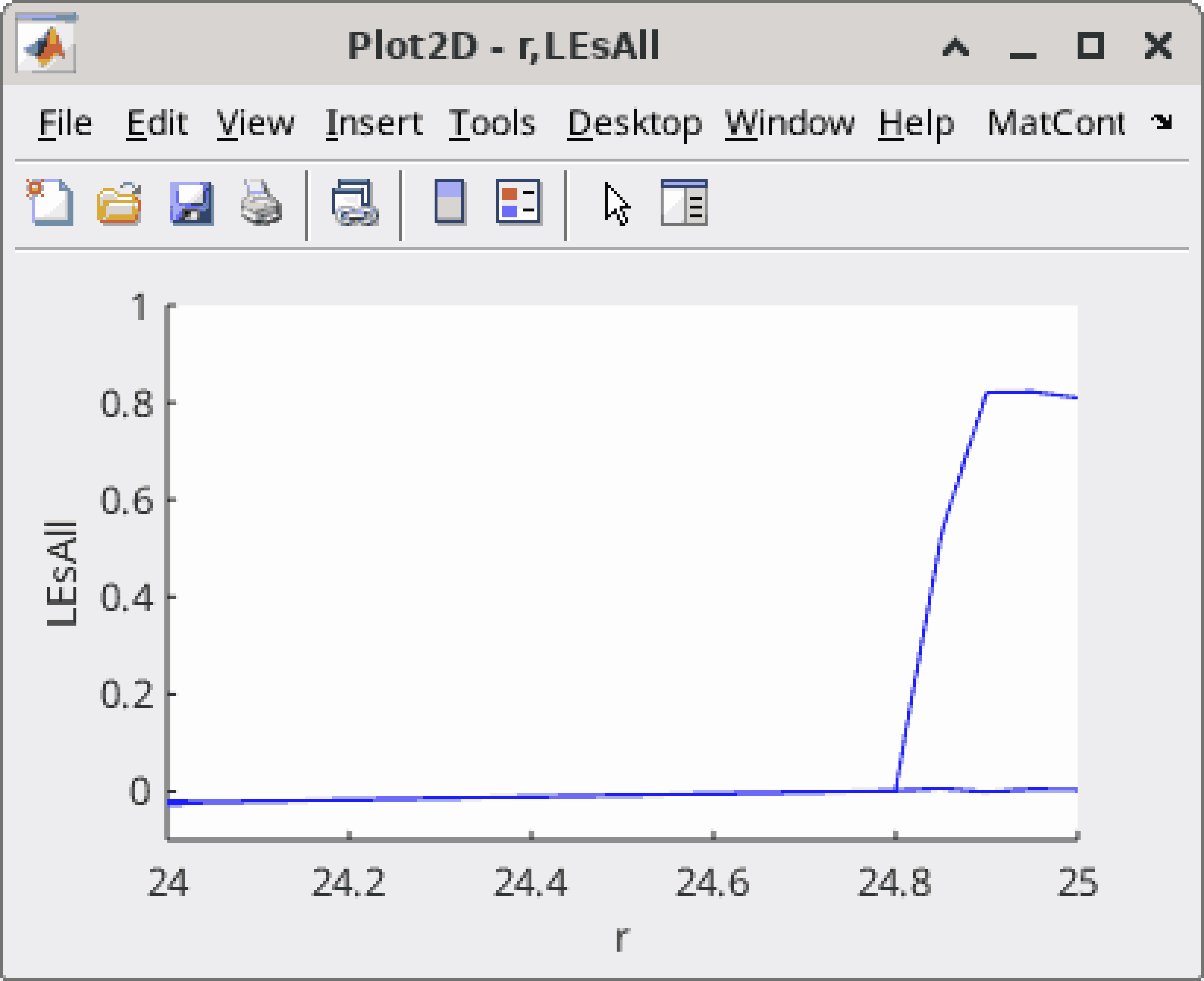}%
    \hfill%
    \includegraphics[width=0.32\linewidth]{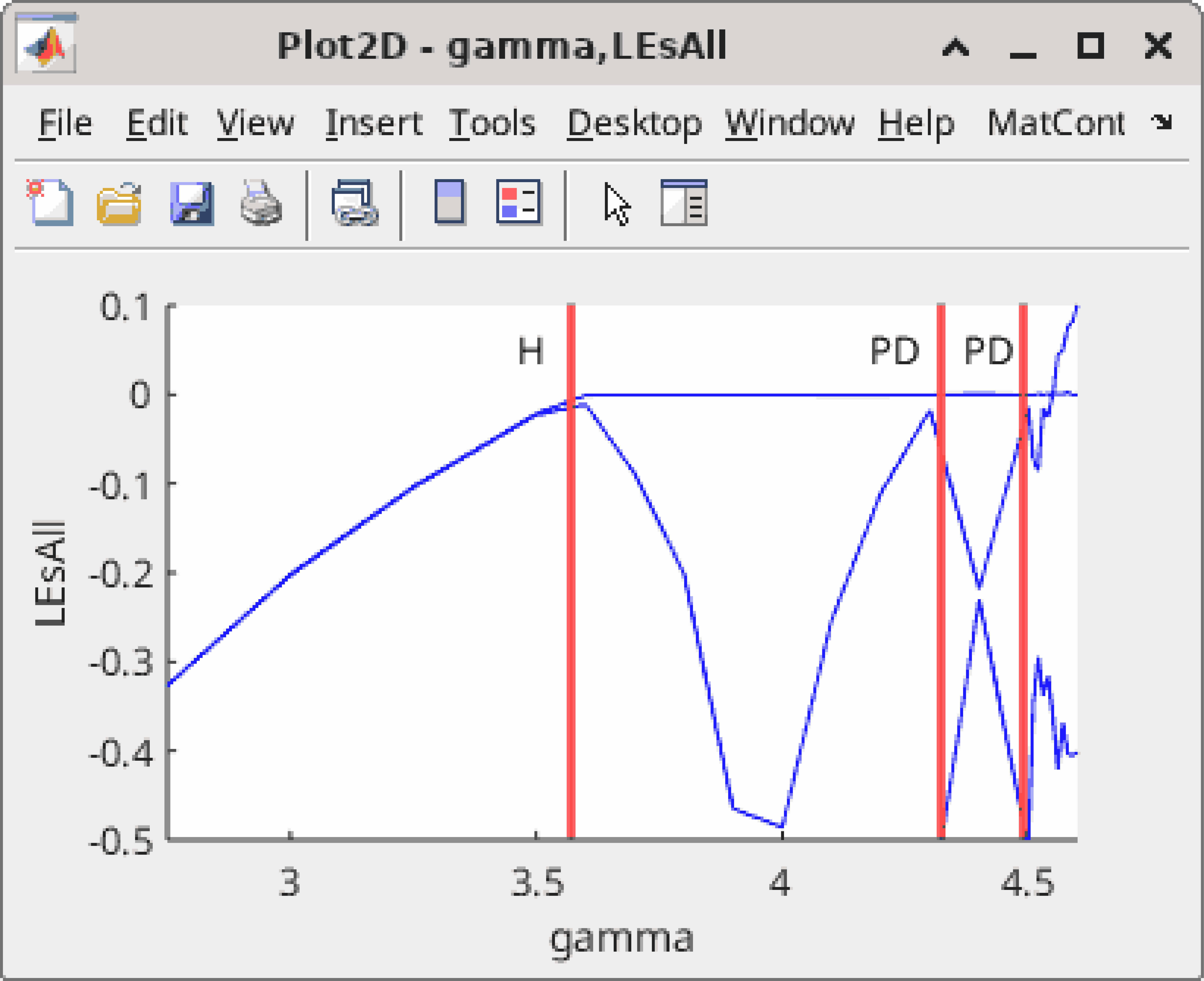}
    \Description{Screenshots of MatCont's \window{Plot2D} windows containing the plots described in the caption.}
    \caption{%
    LEs computed within MatCont.
    Lorenz equations: running estimates of all LEs showing convergence (left), and the two most positive LEs varying the parameter $r$ near the subcritical Hopf bifurcation at $r_{H}\approx 24.74$ (middle), where the jump of LE$_{1}$ indicates the appearance of chaos.
    Quadratic renewal equation \eqref{eq:re-quad}: the two/three most positive LEs using the approximating ODE with $M=10$ (right), varying the parameter $\gamma$ through a period-doubling cascade. Red lines indicate some of the known bifurcations.
    }
    \label{fig:LE-demos}
\end{figure}

To demonstrate the computation of LEs, we turn to the classical Lorenz equations and to a renewal equation. 
For the Lorenz equations, $x'=s(y-x)$, $y'=x(r-z)-y$, $z'=xy-bz$ we set the classical parameters $s=10$ and $b=8/3$, and we use $r$ as a varying parameter.
First, this result shows that the user can inspect the convergence of the computation as the interval is long enough (Figure~\ref{fig:LE-demos}, left, $r=28$). Second, we can see the evolution of the exponents as we follow an attractor and observe bifurcations (Figure~\ref{fig:LE-demos}, middle, $24\leq r\leq 25$). Third, the renewal equation (described below) demonstrates that we can study this system combining both new additions to MatCont (Figure~\ref{fig:LE-demos}, right).

\paragraph{A renewal equation with quadratic nonlinearity}
The equation
\begin{equation}\label{eq:re-quad}
x(t) = \frac{\gamma}{2}\int_{-3}^{-1} x(t+\theta) (1-x(t+\theta)) d\theta
\end{equation}
is input in the delay equation importer as
\begin{lstlisting}
x[t]=gamma/2*DE_int(@(theta)x[t+theta]*(1-x[t+theta]),-3,-1)
\end{lstlisting}
with \lstinline!x! as the list of coordinates and \lstinline!gamma! as the list of parameters.
It exhibits a transcritical bifurcation of the trivial and nontrivial equilibria at $\gamma=1$, and a Hopf bifurcation of the nontrivial equilibrium at $\gamma=2+\frac{\pi}{2}$ from which a branch of stable limit cycles of period $4$ arises.
Experimentally, it is known to also have subsequent period-doubling bifurcations at $\gamma\approx 4.32,4.49,4.53$, and it appears to enter a chaotic regime for $\gamma\geq4.55$. 
See \cite{BredaDiekmannLiessiScarabel2016,BredaLiessi2018,BredaLiessi2024} for more details and other properties of~\eqref{eq:re-quad}.
We have computed the LEs within MatCont: Figure \ref{fig:LE-demos} (right) shows the successive bifurcations.
The computation time, with $M=10$ and truncation time $1000$ for the QR method, was $813.28$\,s with MATLAB R2024b on an HP Probook 430 G7 with Intel Core i5-10210U CPU ($1.60$\,GHz) and $16$\,GB of RAM running Fedora Linux 40. %DL's machine
% Lorenz1 2.89 s on DL's machine
% Lorenz2 56.31 s on DL's machine

\begin{figure}
    \centering
    \includegraphics[width=0.45\linewidth]{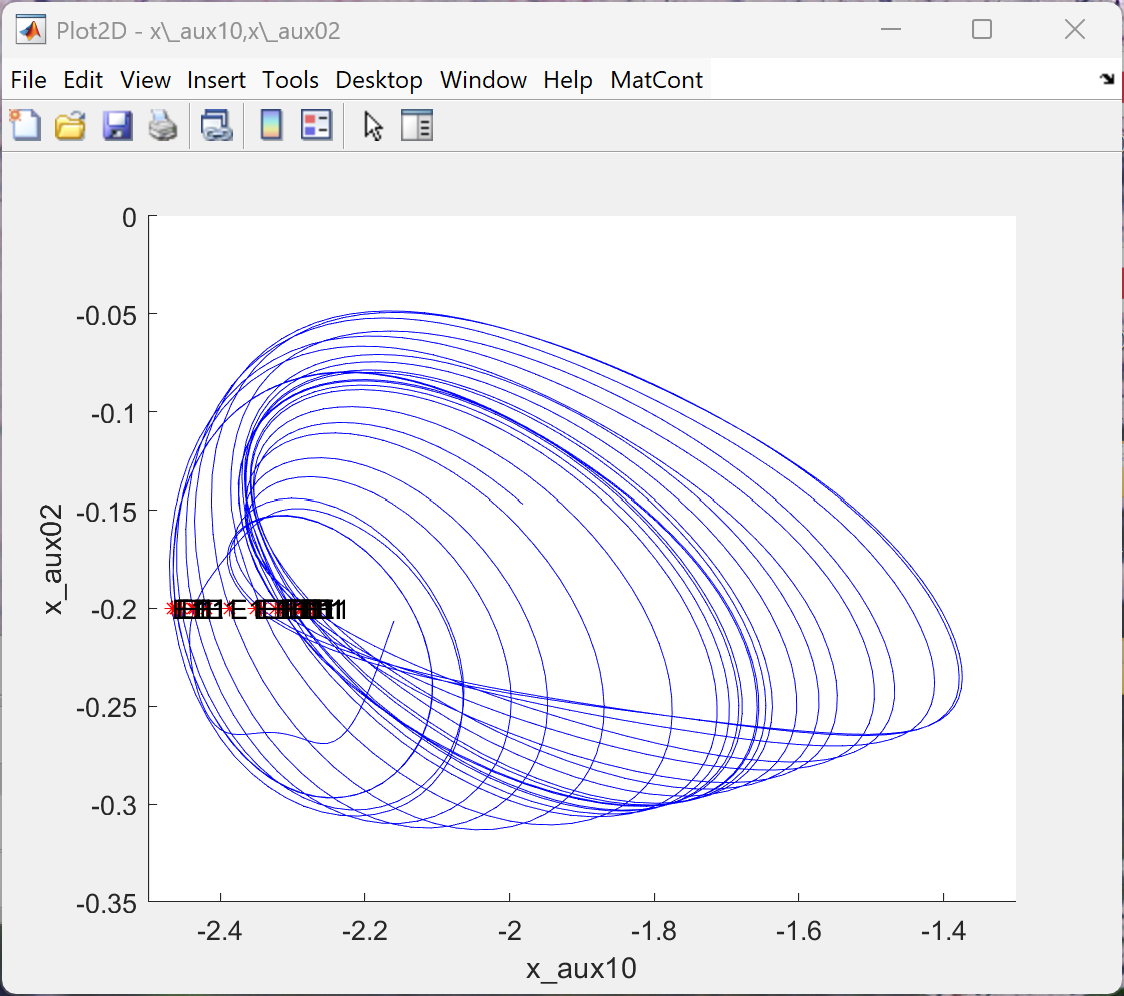}%
    \hfill%
    \includegraphics[width=0.45\linewidth]{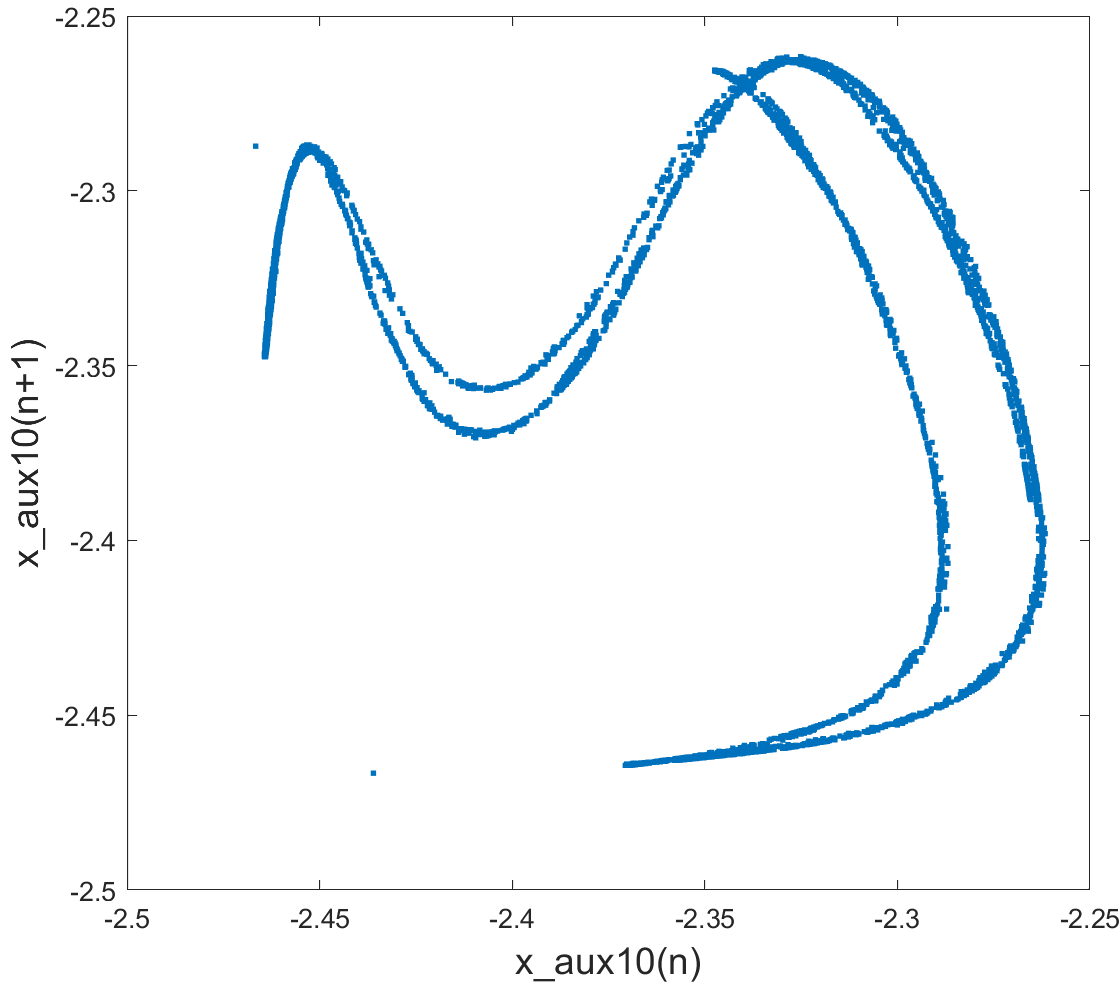}%
    \Description{Screenshots of MatCont's \window{Plot2D} windows containing the plots described in the caption.}
    \caption{%
     Left: Projection of the orbit onto the $(\text{\lstinline!x_aux10!},\text{\lstinline!x_aux02!})$-plane with $t_{\text{max}}=100$. Every crossing is indicated by a marker and label \guielement{E1}. Right: Recurrence plot for the first component \lstinline!x_aux10! for a trajectory with $t_{\text{max}}=20000$. The bending and self-similarity are typical of chaos.
    }
    \label{fig:Event-demos}
\end{figure}

Next, we want to illustrate for $\gamma=4.6$ that the dynamics on the Poincar\'e section exhibit properties typical of chaotic dynamics. We choose the section $\text{\lstinline!x_aux02!}+0.2=0$ and in the window \window{Integrator} we enter \lstinline!MyEvent! to the field \field{Eventfunction} which refers to an m-file with the following code. We specify that we are only interested in crossings in one direction. In the 2D plot, the section appears as desired, see Figure \ref{fig:Event-demos} (left panel). Extending the orbit, and taking only the \lstinline!x_aux10! coordinate on the section, we can create a recurrence plot, see Figure \ref{fig:Event-demos} (right panel), which illustrates the chaotic dynamics expected from a positive Lyapunov exponent. For the purpose of illustrating the dynamics, we show the auxiliary variables, rather the reconstructed one needed for model interpretation.

\begin{lstlisting}
function [res,isterminal,direction]=MyEvent(t,x,varargin)
res=x(2)+0.2; %Definining the section
isterminal=0; %Do not stop
direction=1;  %0 for detecting all crossings, +1/-1 for in/decreasing
\end{lstlisting}

%for RE, we used g=[2.75:.25:3.5 3.6:.1:4.4 4.5:.01:4.6]
% starting from point setup with constant function 0.1
% markers added with the following commands
%hold on
%xline(3.5707963267949,'Color',[1 0 0],'LineWidth',2);
%xline(4.32,'Color',[1 0 0],'LineWidth',2);
%xline(4.49,'Color',[1 0 0],'LineWidth',2);
%annotation('textbox',[0.41 0.8 0.06 0.12],'String',{'H'},'FitBoxToText','off','EdgeColor','none');
%annotation('textbox',[0.71 0.8 0.06 0.12],'String',{'PD'},'FitBoxToText','off','EdgeColor','none');
%annotation('textbox',[0.79 0.8 0.06 0.12],'String',{'PD'},'FitBoxToText','off','EdgeColor','none');

\subsection{Other examples}
\begin{figure}
    \centering
    \includegraphics[width=0.32\linewidth]{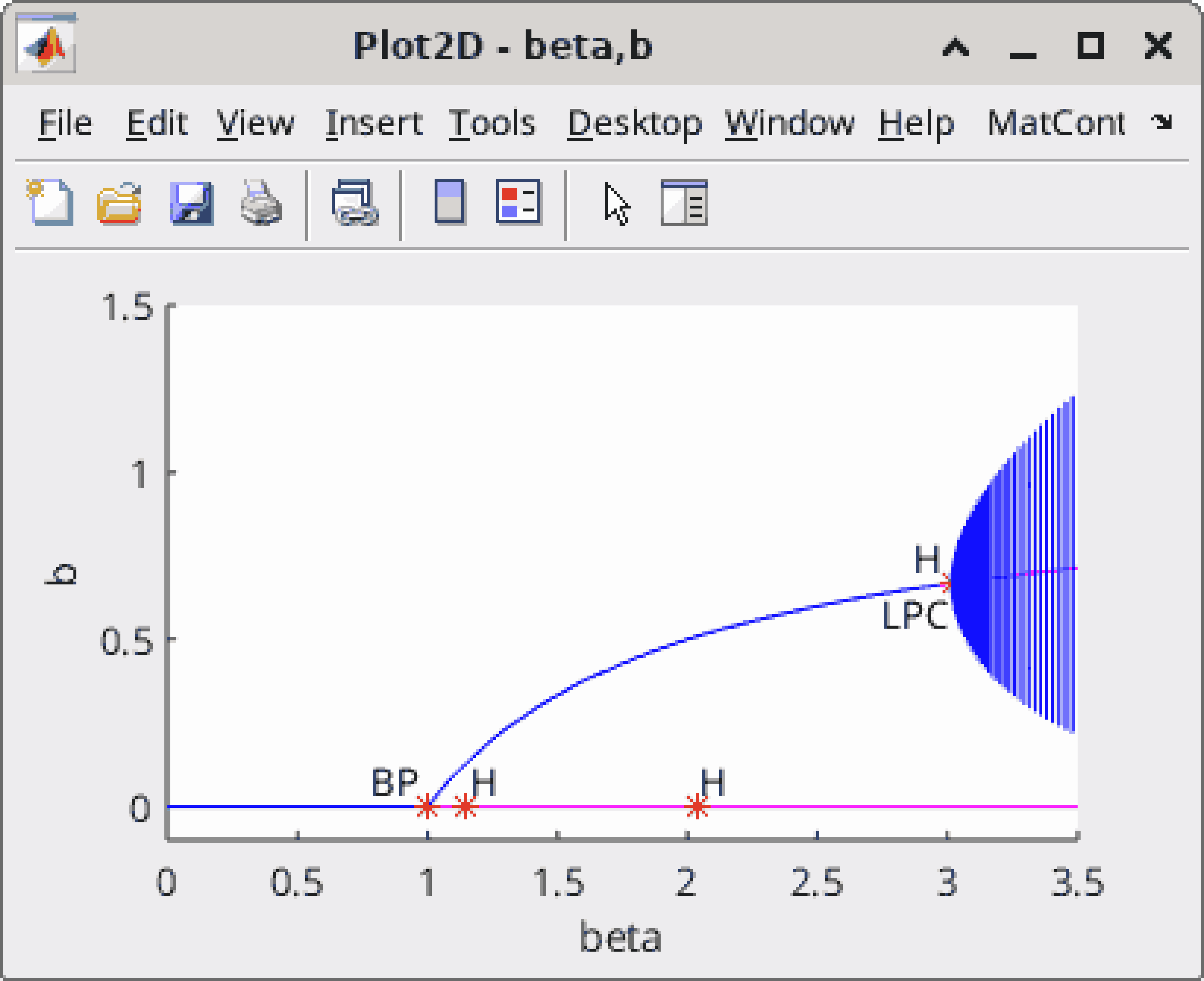}%
    \hfill%
    \includegraphics[width=0.32\linewidth]{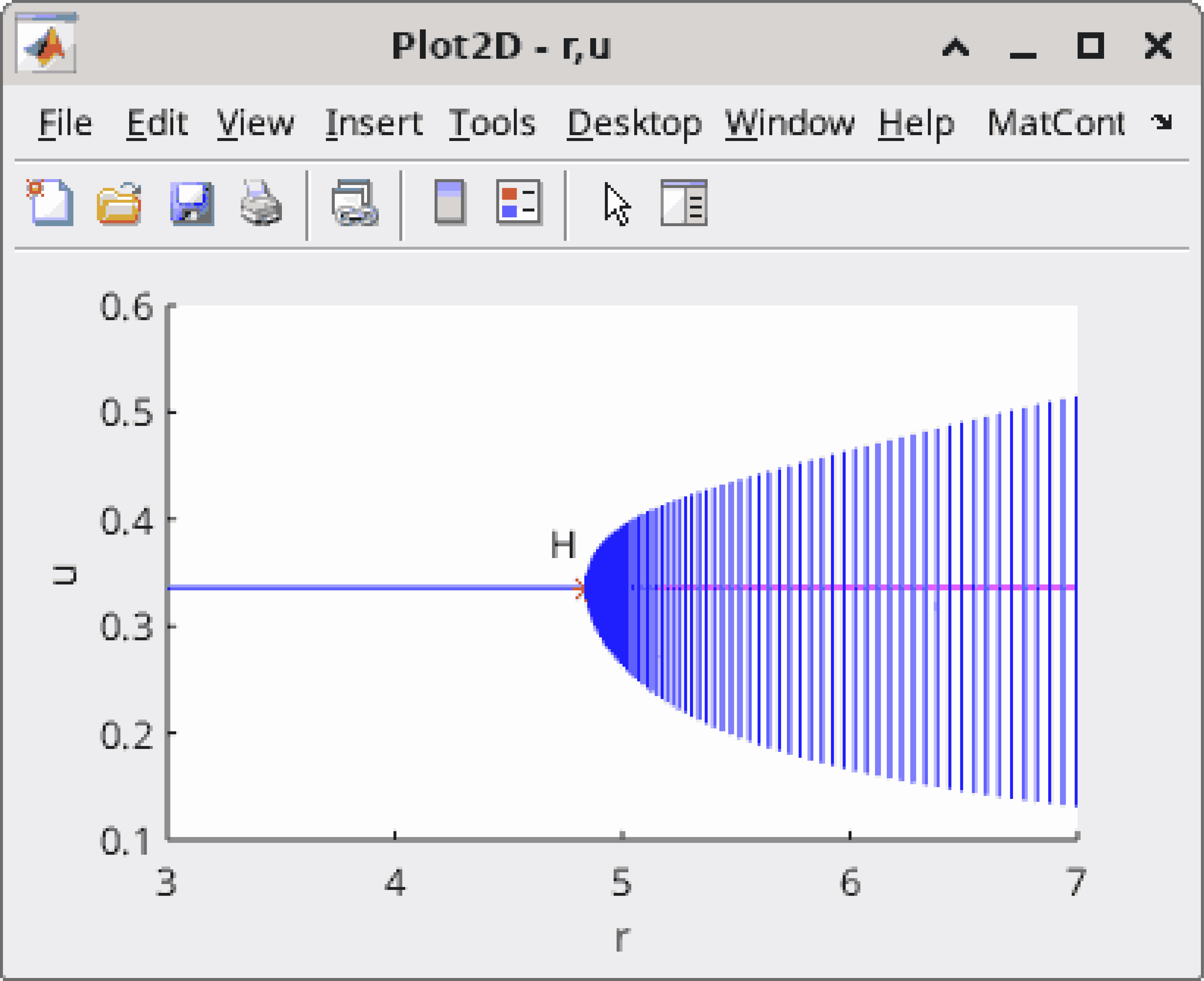}%
    \hfill%
    \includegraphics[width=0.32\linewidth]{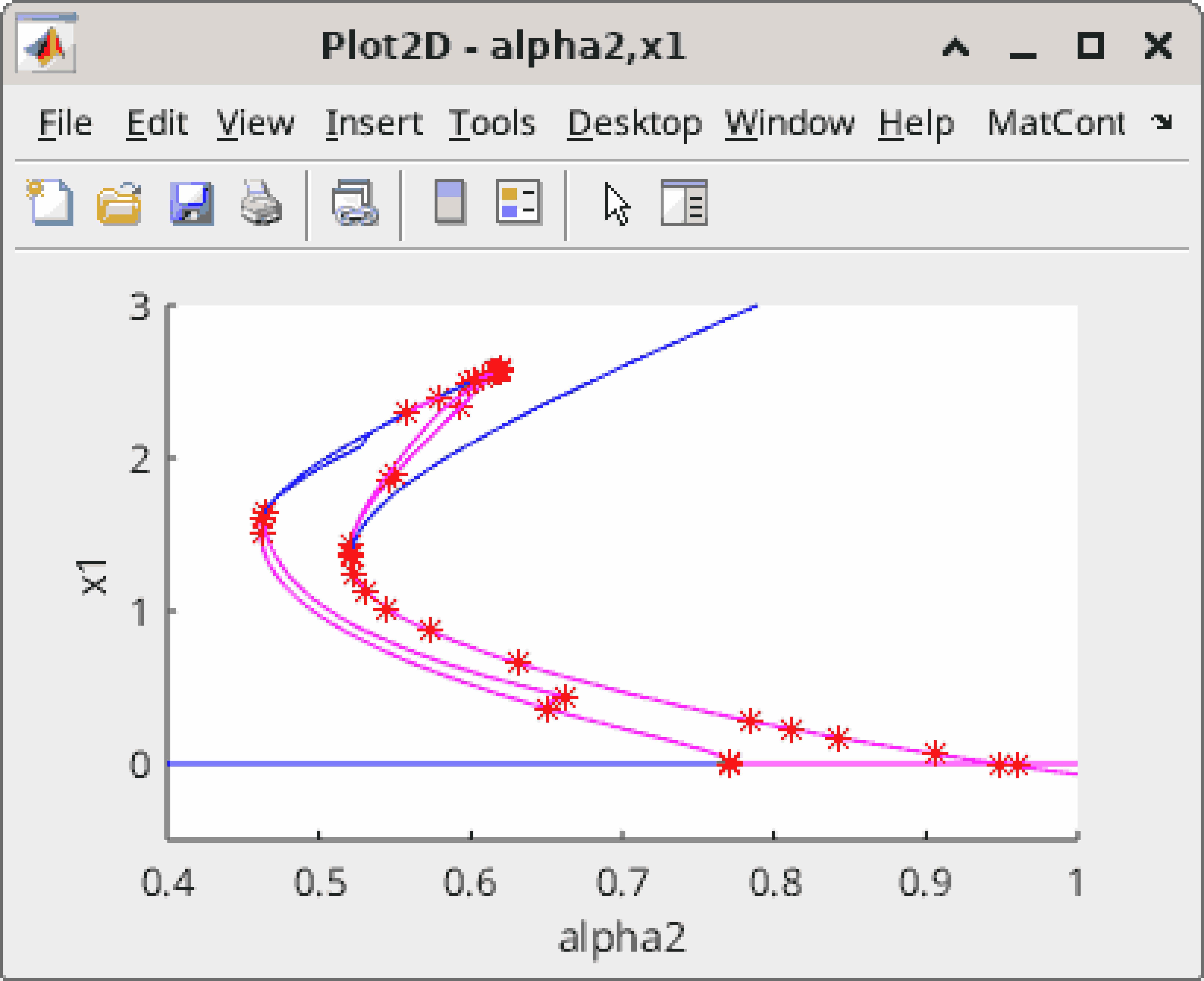}
    \Description{Screenshots of MatCont's \window{Plot2D} windows containing the plots described in the caption.}
    \caption{Bifurcation diagrams for the simplified logistic Daphnia model \eqref{eq:logdaphnia} (left), the refractory neural network model \eqref{eq:refractory} (middle), and the two-node neural network \eqref{eq:neuron2} (right; we removed the point labels for clarity).}
    \label{fig:other-examples}
\end{figure}

As further examples, we collect some DEs, the corresponding input for the delay equation importer, and some known properties that the user can try to reproduce with MatCont.
Figure~\ref{fig:other-examples} shows the bifurcation diagrams obtained with MatCont for these examples.

\paragraph{A simplified logistic Daphnia model}
The system of coupled RE and DDE 
\begin{equation}\label{eq:logdaphnia}
\left\{
\begin{aligned}
& b(t) = \beta S(t) \int_{a_{\text{repr}}}^{a_{\text{max}}} b(t-a) \D a \\
& S'(t) = r S(t) \left(1 - \frac{S(t)}{K}\right) - \gamma S(t) \int_{a_{\text{repr}}}^{a_{\text{max}}} b(t-a) \D a
\end{aligned}
\right.
\end{equation}
is input as
\begin{lstlisting}
S_int_b=S[t]*DE_int(@(a)b[t-a],a_repr,a_max)
b[t]=beta*S_int_b
S'[t]=r*S[t]*(1-S[t]/K)-gamma*S_int_b
\end{lstlisting}
with \lstinline!b,S! as the list of coordinates and \lstinline!a_repr,a_max,beta,r,K,gamma! as the list of parameters.
Note that we stored the value of the repeated integral in a variable for efficiency.
For fixed $a_\text{repr}=3$, $a_{\text{max}}=4$ and $r=K=\gamma=1$, the consumer-free equilibrium $(0,K)=(0,1)$ undergoes a transcritical bifurcation at $\beta=(K(a_{\text{max}}-a_{\text{repr}}))^{-1}=1$, where the nontrivial equilibrium becomes positive and stable.
The latter then exhibits a Hopf bifurcation at $\beta\approx3.0161$.
For more details, see~\cite{BredaDiekmannGyllenbergScarabelVermiglio2016,BredaLiessi2020}.

\paragraph{A refractory neural network model}
In order to input the DDE
\begin{equation*}
u'(t) = r\left(- u(t) + \left( 1 - \int_{t-1}^t u(s) ds \right) \frac{1}{1+e^{-a(u(t)+\theta)}}\right),
\end{equation*}
we need to change variable in the integral to make the integration endpoints independent of $t$:
\begin{equation}\label{eq:refractory}
u'(t) = r\left(- u(t) + \left( 1 - \int_{-1}^0 u(t+s) ds \right) \frac{1}{1+e^{-a(u(t)+\theta)}}\right).
\end{equation}
We can now input
\begin{lstlisting}
u'[t]=r*(-u[t]+(1-DE_int(@(s)u[t+s],-1,0))*1/(1+exp(-a*(u[t]+theta))))
\end{lstlisting}
with \lstinline!u! as the list of coordinates and \lstinline!r,a,theta! as the list of parameters.
If we fix $a=8$ and $\theta=-0.333$, the equilibrium undergoes a Hopf bifurcation at $r\approx4.8394$.
For more details, see \cite{CurtuErmentrout2001}.

\paragraph{A two-node neural network}
The system of two DDEs with two delays
\begin{equation}\label{eq:neuron2}
\left\{
\begin{aligned}
& x_1'(t) = - x_1(t) - \alpha_1 S(\beta_1 x_1(t-\tau_1)) + \alpha_2 S(\beta_2 x_2(t-\tau_2)) \\
& x_2'(t) = -x_2(t) - \alpha_1 S(\beta_1 x_2(t-\tau_1)) + \alpha_2 S(\beta_2 x_1(t-\tau_2))
\end{aligned}
\right.
\end{equation}
with $S(x) = (\tanh(x-a) + \tanh(a)) \cosh(a)^2$ can be input as
\begin{lstlisting}
S=@(x)(tanh(x-a)+tanh(a))*cosh(a)^2
x1'=-x1[t]-alpha1*S(beta1*x1[t-tau1])+alpha2*S(beta2*x2[t-tau2])
x2'=-x2[t]-alpha1*S(beta1*x2[t-tau1])+alpha2*S(beta2*x1[t-tau2])
\end{lstlisting}
with \lstinline!x1,x2! as the list of coordinates and \lstinline!a,tau1,tau2,alpha1,alpha2,beta1,beta2! as the list of parameters.
Note that, as done for \eqref{eq:logdaphnia}, we stored the definition of the function $S$ for efficiency.
For fixed $a=1$, $\tau_1=11.6$, $\tau_2=20.3$, $\alpha_1=0.069$, $\beta_1=2$, and $\beta_2=1.2$, the trivial equilibrium undergoes a Hopf bifurcation at $\alpha_2\approx0.771$ and a transcritical bifurcation at $\alpha_2\approx0.948$.
From there, using collocation polynomials of degree $20$, a rather rich bifurcation diagram can be computed; see~\cite{VisserMeijerVanPuttenVanGils2012} for more details.

\section{Concluding remarks}

The new release 7p6 extends the capabilities of MatCont allowing users to:
\begin{itemize}
    \item analyse DEs with finite delay, through a delay equation importer that automatically implements the pseudospectral approximation via an ODE; 
    \item compute LEs. 
\end{itemize}

In this paper we have included some examples that illustrate the new functionalities, and we refer the interested user to \url{https://sourceforge.net/projects/matcont/} for extended tutorials. 
We have also included a performance comparison on the Mackey--Glass equation between the MatCont implementation and DDE-BIFTOOL. It showed that analysing DEs with MatCont can become computationally expensive, since the underlying ODE approximation has a large dimension. However, the pseudospectral method allows to obtain approximation with reasonable accuracy for low orders compared to other approximation methods such as the linear chain trick. 
Moreover, some types of equations, like integral equations, are currently supported by the MatCont delay equations importer but not by DDE-BIFTOOL or other packages. 

The exponential convergence of the stability of equilibria, for DDEs and REs, has been proved in \cite{BredaMasetVermiglio2015,ScarabelDiekmannVermiglio2021}. MatCont also allows to perform time integration and analyse more complex bifurcations, but it is worth stressing that the theoretical convergence of these operations is still an open problem. Users should validate the numerical results by testing approximations of increasing orders, as we demonstrated in Section \ref{sec:mgconvergence}. 

Some limitations of the delay equation importer could be improved in subsequent versions. In particular, user functions are currently not supported for DEs. Moreover, ODEs in a delay system are discretised like DDEs introducing additional auxiliary variables, although in principle there is no need to do so. Additional functionalities we aim to add are the support of infinite delays, following ideas presented in \cite{GyllenbergScarabelVermiglio2018, ScarabelVermiglio2024}, and state-dependent delays, for which a preliminary investigation has been conducted in \cite{GettoGyllenbergNakataScarabel2019}, where the delay can be implicitly defined via a condition that involves the integral of the state.
We also aim to allow the personalised selection of different discretisation grids.  

We also note that importers are a flexible and general solution in the long term: they can be maintained and extended independently of MatCont proper and new importers can easily be introduced to discretise other kinds of abstract differential equations into ODEs (including, e.g., PDEs with vector variables due to discretisation).
A similar approach was used in previous versions of MatCont for importing SBML%
\footnote{Systems Biology Markup Language, \url{http://sbml.org/}} files.

Finally, we stress that a large collection of other software packages exists for continuation and bifurcation analyses of ODEs, e.g., AUTO \cite{AUTO} (and XPP), and PyDSTool \cite{PyDSTool}, to name a few. 
By reducing the DE to an ODE, the pseudospectral approximation is in principle compatible with any of these packages. The advantage of the current MatCont importer is that it offers an easy-to-use GUI for the analyses. 
In a similar spirit, we implemented the computation of LEs. One can find open-source codes for these computations, but using them requires another implementation of the system, which one preferably avoids. Again, the MatCont GUI offers a direct approach to this method.

%%%%% FROM ACM TEMPLATE

\begin{acks}
FS and RV are members of INdAM research group GNCS. DL, FS and RV are members of UMI research group `Mo\-del\-li\-sti\-ca socio-epidemiologica'. FS is a member of JUNIPER (Joint UNIversities Pandemic and Epidemiological Research).
The project of the delay equation importer was started as part of ES's BSc internship and thesis \cite{Santi2022} at CDLab supervised by RV and DL. MT developed code for LE computations supervised by HM.
Most of this research has been conducted during research visits at the University of Leeds, funded by a Small Grant from the Heilbronn Institute for Mathematical Research (HIMR) and by the University of Leeds Research Equity, Diversity and Inclusion (REDI) Pilot Fund, as well as during the workshop `Towards rigorous results in state-dependent delay equations' held at the Lorentz Center in Leiden on 4--8 March 2024.
The authors also received partial financial support
by GNCS projects 2023 `Sistemi dinamici e modelli di evoluzione: tecniche funzionali, analisi qualitativa e metodi numerici'
(CUP: E53C22001930001)
and 2024 `Analisi numerica di problemi di evoluzione complessi: stabilità, conservazione e tecniche data-driven'
(CUP: E53C23001670001),
Italian Ministry of University and Research (MUR)
project PRIN 2020 (No.\ 2020JLWP23) `Integrated Mathematical Approaches to Socio-Epidemiological Dynamics', Unit of Udine (CUP: G25F22000430006),
QJMAM Fund for Applied Mathematics,
European Union NextGenerationEU project `MOdellistica Numerica e Data-driven per l'Innovazione sostenibile -- MONDI'
CUP: G25F21003390007,
Spanish AEI project `Modelización y análisis numérico en problemas de evolución con aplicaciones a biología, economía y mecánica de fluídos'
PID2020-113554GB-I00/AEI/10.13039/501100011033.

Author contributions:
conceptualisation, methodology, writing -- review \& editing: DL, HM, FS, RV;
funding acquisition: FS, RV;
software: DL, HM, ES, MT;
writing -- original draft: DL.
%% https://credit.niso.org/
HM and FS are joint last authors.
\end{acks}

\printbibliography

@Article{AndoBredaGava2020,
  author       = {Andò, Alessia and Breda, Dimitri and Gava, Giulia},
  title        = {How fast is the linear chain trick? A rigorous analysis in the context of behavioral epidemiology},
  journaltitle = {Mathematical Biosciences and Engineering},
  shortjournal = {Math. Biosci. Eng.},
  volume       = {17},
  number       = {5},
  date         = {2020-07-24},
  pages        = {5059-5084},
  doi          = {10.3934/mbe.2020273},
}

@Software{AUTO,
  author = {Doedel, Eusebius J. and Oldeman, Bart E.},
  title  = {{AUTO-07P}: Continuation and bifurcation software for ordinary differential equations},
  date   = {2021-08-25},
  url    = {http://indy.cs.concordia.ca/auto/},
}

@Article{BerrutTrefethen2004,
  author       = {Berrut, Jean-Paul and Trefethen, Lloyd Nicholas},
  title        = {Barycentric {L}agrange Interpolation},
  journaltitle = {SIAM Review},
  shortjournal = {SIAM Rev.},
  date         = {2004},
  volume       = {46},
  number       = {3},
  pages        = {501-517},
  doi          = {10.1137/S0036144502417715},
}

@Article{Bosschaert:2024,
  author       = {Bosschaert, Maikel M. and Kuznetsov, Yuri A.},
  title        = {Interplay between Normal Forms and Center Manifold Reduction for Homoclinic Predictors near {B}ogdanov–{T}akens Bifurcation},
  journaltitle = {SIAM Journal on Applied Dynamical Systems},
  shortjournal = {SIAM J. Appl. Dyn. Syst.},
  volume       = {23},
  number       = {1},
  pages        = {410-439},
  date         = {2024},
  doi          = {10.1137/22M151354X},
}

@Article{BredaDiekmannGyllenbergScarabelVermiglio2016,
  author       = {Breda, Dimitri and Diekmann, Odo and Gyllenberg, Mats and Scarabel, Francesca and Vermiglio, Rossana},
  title        = {Pseudospectral discretization of nonlinear delay equations: New prospects for numerical bifurcation analysis},
  journaltitle = {SIAM Journal on Applied Dynamical Systems},
  shortjournal = {SIAM J. Appl. Dyn. Syst.},
  volume       = {15},
  number       = {1},
  date         = {2016},
  pages        = {1-23},
  doi          = {10.1137/15M1040931},
}

@Article{BredaLiessi2018,
  author       = {Breda, Dimitri and Liessi, Davide},
  title        = {Approximation of eigenvalues of evolution operators for linear renewal equations},
  journaltitle = {SIAM Journal on Numerical Analysis},
  shortjournal = {SIAM J. Numer. Anal.},
  date         = {2018},
  volume       = {56},
  number       = {3},
  pages        = {1456-1481},
  doi          = {10.1137/17M1140534},
}

@Article{BredaLiessi2020,
  author       = {Breda, Dimitri and Liessi, Davide},
  title        = {Approximation of eigenvalues of evolution operators for linear coupled renewal and retarded functional differential equations},
  journaltitle = {Ricerche di Matematica},
  shortjournal = {Ric. Mat.},
  date         = {2020-11},
  volume       = {69},
  number       = {2},
  pages        = {457-481},
  doi          = {10.1007/s11587-020-00513-9},
}

@Article{BredaLiessi2021,
  author       = {Breda, Dimitri and Liessi, Davide},
  title        = {Floquet theory and stability of periodic solutions of renewal equations},
  journaltitle = {Journal of Dynamics and Differential Equations},
  shortjournal = {J. Dynam. Differential Equations},
  volume       = {33},
  number       = {2},
  date         = {2021-06},
  pages        = {677-714},
  doi          = {10.1007/s10884-020-09826-7},
}

@Article{BredaLiessi2024,
  author       = {Breda, Dimitri and Liessi, Davide},
  title        = {A practical approach to computing {L}yapunov exponents of renewal and delay equations},
  journaltitle = {Mathematical Biosciences and Engineering},
  shortjournal = {Math. Biosci. Eng.},
  volume       = {21},
  number       = {1},
  pages        = {1249-1269},
  date         = {2024},
  doi          = {10.3934/mbe.2024053},
}

@Book{BredaMasetVermiglio2015,
  author      = {Breda, Dimitri and Maset, Stefano and Vermiglio, Rossana},
  title       = {Stability of Linear Delay Differential Equations},
  date        = {2015},
  shortseries = {SpringerBriefs Control Autom. Robot.},
  series      = {SpringerBriefs in Control, Automation and Robotics},
  subtitle    = {A Numerical Approach with MATLAB},
  publisher   = {Springer},
  location    = {New York},
  doi         = {10.1007/978-1-4939-2107-2},
}

@Article{BredaDiekmannLiessiScarabel2016,
  author       = {Breda, Dimitri and Diekmann, Odo and Liessi, Davide and Scarabel, Francesca},
  title        = {Numerical bifurcation analysis of a class of nonlinear renewal equations},
  journaltitle = {Electronic Journal of Qualitative Theory of Differential Equations},
  shortjournal = {Electron. J. Qual. Theory Differ. Equ.},
  volume       = {2016},
  eid          = {65},
  date         = {2016-09-12},
  pages        = {1-24},
  doi          = {10.14232/ejqtde.2016.1.65},
}

@Article{CurtuErmentrout2001,
  author       = {Curtu, Rodica and Ermentrout, Bard},
  title        = {Oscillations in a refractory neural net},
  journaltitle = {Journal of Mathematical Biology},
  shortjournal = {J. Math. Biol.},
  volume       = {43},
  pages        = {81-100},
  date         = {2001-07},
  doi          = {10.1007/s002850100089},
}

@Article{ClenshawCurtis1960,
  author       = {Clenshaw, C. W. and Curtis, A. R.},
  title        = {A method for numerical integration on an automatic computer},
  journaltitle = {Numerische Mathematik},
  shortjournal = {SIAM J. Numer. Anal.},
  volume       = {2},
  date         = {1960-12},
  pages        = {197-205},
  doi          = {10.1007/BF01386223},
}

@Book{Cushing1977,
  author      = {Cushing, Jim M.},
  title       = {Integrodifferential Equations and Delay Models in Population Dynamics},
  date        = {1977},
  shortseries = {Lect. Notes Biomath.},
  series      = {Lecture Notes in Biomathematics},
  number      = {20},
  publisher   = {Springer},
  location    = {Berlin, Heidelberg},
  doi         = {10.1007/978-3-642-93073-7},
}

@Article{DeWolffScarabelVerduynLunelDiekmann2021,
  author       = {de Wolff, Babette and Scarabel, Francesca and Verduyn Lunel, Sjoerd and Diekmann, Odo},
  title        = {Pseudospectral approximation of {H}opf bifurcation for delay differential equations},
  journaltitle = {SIAM Journal on Applied Dynamical Systems},
  shortjournal = {SIAM J. Appl. Dyn. Syst.},
  volume       = {20},
  number       = {1},
  date         = {2021-03-01},
  pages        = {333-370},
  doi          = {10.1137/20M1347577},
  eprinttype   = {arXiv},
  eprintclass  = {math.DS},
  eprint       = {2006.13810},
}

@Article{DhoogeGovaertsKuznetsov2003,
  author       = {Dhooge, A. and Govaerts, W. and Kuznetsov, Yuri A.},
  title        = {{MATCONT}: {A} {MATLAB} package for numerical bifurcation analysis of {ODEs}},
  journaltitle = {ACM Transactions on Mathematical Software},
  shortjournal = {ACM Trans. Math. Software},
  date         = {2003-06},
  volume       = {29},
  number       = {2},
  pages        = {141-164},
  doi          = {10.1145/779359.779362},
}

@Article{DiekmannGettoGyllenberg2008,
  author       = {Diekmann, Odo and Getto, Philipp and Gyllenberg, Mats},
  title        = {Stability and bifurcation analysis of {V}olterra functional equations in the light of suns and stars},
  journaltitle = {SIAM Journal on Mathematical Analysis},
  shortjournal = {SIAM J. Math. Anal.},
  date         = {2008},
  volume       = {39},
  number       = {4},
  pages        = {1023-1069},
  doi          = {10.1137/060659211},
}

@Article{DiekmannScarabelVermiglio2020,
  author       = {Diekmann, Odo and Scarabel, Francesca and Vermiglio, Rossana},
  title        = {Pseudospectral discretization of delay differential equations in sun-star formulation: Results and conjectures},
  journaltitle = {Discrete and Continuous Dynamical Systems. Series S},
  shortjournal = {Discrete Contin. Dyn. Syst. Ser. S},
  volume       = {13},
  number       = {9},
  date         = {2020-09},
  pages        = {2575-2602},
  doi          = {10.3934/dcdss.2020196},
}

@Book{DiekmannVanGilsVerduynLunelWalther1995,
  author      = {Diekmann, Odo and van Gils, Stephan A. and Verduyn Lunel, Sjoerd M. and Walther, Hans-Otto},
  title       = {Delay Equations},
  date        = {1995},
  shortseries = {Appl. Math. Sci.},
  number      = {110},
  series      = {Applied Mathematical Sciences},
  subtitle    = {Functional-, Complex- and Nonlinear Analysis},
  publisher   = {Springer},
  location    = {New York},
  doi         = {10.1007/978-1-4612-4206-2},
}

@Article{EngelborghsLuzyaninaRoose2002,
  author       = {Engelborghs, Koen and Luzyanina, Tatyana and Roose, Dirk},
  title        = {Numerical Bifurcation Analysis of Delay Differential Equations Using {DDE-BIFTOOL}},
  journaltitle = {ACM Transactions on Mathematical Software},
  shortjournal = {ACM Trans. Math. Softw.},
  date         = {2002},
  volume       = {28},
  number       = {1},
  pages        = {1--21},
  doi          = {10.1145/513001.513002},
}

@Book{Ermentrout2002,
  author      = {Ermentrout, Bard},
  title       = {Simulating, Analyzing, and Animating Dynamical Systems: A Guide to {XPPAUT} for Researchers and Students},
  date        = {2002},
  shortseries = {Software Environ. Tools},
  series      = {Software, Environments and Tools},
  number      = {14},
  publisher   = {Society for Industrial and Applied Mathematics},
  location    = {Philadelphia},
  doi         = {10.1137/1.9780898718195},
}

@Article{GettoGyllenbergNakataScarabel2019,
  author       = {Getto, Philipp and Gyllenberg, Mats and Nakata, Yukihiko and Scarabel, Francesca},
  title        = {Stability analysis of a state-dependent delay differential equation for cell maturation: Analytical and numerical methods},
  journaltitle = {Journal of Mathematical Biology},
  shortjournal = {J. Math. Biol.},
  volume       = {79},
  number       = {1},
  date         = {2019-07},
  pages        = {281-328},
  doi          = {10.1007/s00285-019-01357-0},
}

@Software{Govorukhin:2004,
  author = {Govorukhin, Vasiliy},
  title  = {Calculation Lyapunov Exponents for ODE},
  date   = {2004-03-18},
  organization   = {MATLAB Central File Exchange},
  url    = {https://www.mathworks.com/matlabcentral/fileexchange/4628-calculation-lyapunov-exponents-for-ode},
}

@Article{GyllenbergScarabelVermiglio2018,
  author       = {Gyllenberg, Mats and Scarabel, Francesca and Vermiglio, Rossana},
  title        = {Equations with infinite delay: Numerical bifurcation analysis via pseudospectral discretization},
  journaltitle = {Applied Mathematics and Computation},
  shortjournal = {Appl. Math. Comput.},
  volume       = {333},
  date         = {2018-09-15},
  pages        = {490-505},
  doi          = {10.1016/j.amc.2018.03.104},
}

@Article{FFJin:1997,
  author       = {Jin, Fei-Fei},
  title        = {An equatorial ocean recharge paradigm for {ENSO}. {P}art {I}: {C}onceptual model},
  journaltitle = {Journal of the Atmospheric Sciences},
  shortjournal = {J. Atmos. Sci.},
  date         = {1997-04-01},
  pages        = {811-829},
  doi          = {10.1175/1520-0469(1997)054<0811:AEORPF>2.0.CO;2},
}

@Incollection{KrauskopfSieber2023,
  author       = {Krauskopf, Bernd and Sieber, Jan},
  title        = {Bifurcation Analysis of Systems With Delays: Methods and Their Use in Applications},
  booktitle    = {Controlling Delayed Dynamics},
  booksubtitle = {Advances in Theory, Methods and Applications},
  editor       = {Breda, Dimitri},
  series       = {CISM International Centre for Mechanical Sciences},
  shortseries  = {CISM},
  number       = {604},
  publisher    = {Springer},
  location     = {Cham},
  date         = {2023},
  pages        = {195-245},
  doi          = {10.1007/978-3-031-01129-0_7},
}

@Book{Kuznetsov:2023,
  author      = {Kuznetsov, Yuri A.},
  title       = {Elements of Applied Bifurcation Theory},
  date        = {2023},
  edition     = {4},
  shortseries = {Appl. Math. Sci.},
  number      = {112},
  series      = {Applied Mathematical Sciences},
  publisher   = {Springer},
  location    = {Cham},
  doi         = {10.1007/978-3-031-22007-4},
}

@Book{KuznetsovMeijer:2019,
  author      = {Kuznetsov, Yuri A. and Meijer, Hil G. E.},
  title       = {Numerical Bifurcation Analysis of Maps},
  subtitle    = {From Theory to Software},
  date        = {2019},
  publisher   = {Cambridge University Press},
  location    = {Cambridge},
  doi         = {10.1017/9781108585804},
}

@Article{MackeyGlass1977,
  author       = {Mackey, Michael C. and Glass, Leon},
  title        = {Oscillation and chaos in physiological control systems},
  journaltitle = {Science},
  shortjournal = {Science},
  volume       = {197},
  number       = {4300},
  pages        = {287-289},
  date         = {1977-07-15},
  doi          = {10.1126/science.267326},
}

@Book{MastroianniMilovanovic2008,
  author      = {Mastroianni, Giuseppe and Milovanović, Gradimir V.},
  title       = {Interpolation Processes},
  date        = {2008-09-25},
  shortseries = {Springer Monogr. Math.},
  series      = {Springer Monographs in Mathematics},
  subtitle    = {Basic Theory and Applications},
  publisher   = {Springer},
  location    = {Berlin, Heidelberg},
  doi         = {10.1007/978-3-540-68349-0},
}

@Article{NevermannGros2023,
  author       = {Nevermann, Daniel Henrik and Gros, Claudius},
  title        = {Mapping dynamical systems with distributed time delays to sets of ordinary differential equations},
  journaltitle = {Journal of Physics A: Mathematical and Theoretical},
  shortjournal = {J. Phys. A},
  volume       = {56},
  number       = {34},
  eid          = {345702},
  date         = {2023-08-03},
  doi          = {10.1088/1751-8121/acea06},
}

@Article{Pusuluri:2021,
  author       = {Pusuluri, K. and Meijer, H. G. E. and Shilnikov, A. L.},
  title        = {Homoclinic puzzles and chaos in a nonlinear laser model},
  journaltitle = {Communications in Nonlinear Science and Numerical Simulation},
  shortjournal = {Commun. Nonlinear Sci. Numer. Simul.},
  volume       = {93},
  eid          = {105503},
  date         = {2021-02},
  doi          = {10.1016/j.cnsns.2020.105503},
}

@Software{PyDSTool,
  author = {Clewley, R. H. and Sherwood, W. E. and LaMar, M. D. and Guckenheimer, J. M.},
  title  = {{PyDSTool}, a software environment for dynamical systems modeling},
  date   = {2007},
  url    = {http://pydstool.sourceforge.net/},
}

@Thesis{Santi2022,
  author      = {Santi, Enrico},
  title       = {Numerical bifurcation analysis of delay equations: a user-friendly extension of {MatCont}'s interface},
  type        = {BSc thesis},
  institution = {University of Udine},
  date        = {2022},
}

@Article{ScarabelDiekmannVermiglio2021,
  author       = {Scarabel, Francesca and Diekmann, Odo and Vermiglio, Rossana},
  title        = {Numerical bifurcation analysis of renewal equations via pseudospectral approximation},
  journaltitle = {Journal of Computational and Applied Mathematics},
  shortjournal = {J. Comput. Appl. Math.},
  volume       = {397},
  eid          = {113611},
  date         = {2021-12-01},
  doi          = {10.1016/j.cam.2021.113611},
  eprinttype   = {arXiv},
  eprintclass  = {math.NA},
  eprint       = {2012.05364},
}

@Article{ScarabelVermiglio2024,
  author       = {Scarabel, Francesca and Vermiglio, Rossana},
  title        = {Equations with infinite delay: pseudospectral discretization for numerical stability and bifurcation in an abstract framework},
  journaltitle = {SIAM Journal on Numerical Analysis},
  shortjournal = {SIAM J. Numer. Anal.},
  volume       = {62},
  number       = {4},
  pages        = {1736-1758},
  date         = {2024-08},
  doi          = {10.1137/23M1581133},
  eprinttype   = {arXiv},
  eprintclass  = {math.NA},
  eprint       = {2306.13351},
}

@article{Scully:2025,
  author       = {Scully, James and Hinsley, Carter and Bloom, David and Meijer, Hil G. E. and Shilnikov, Andrey L.},
  title        = {Widespread neuronal chaos induced by slow oscillating currents},
  journaltitle = {Chaos: An Interdisciplinary Journal of Nonlinear Science},
  shortjournal = {Chaos},
  volume       = {35},
  number       = {3},
  eid          = {033120},
  date         = {2025-03-07},
  doi          = {10.1063/5.0248001},
}

@Book{Smith2011,
  author      = {Smith, Hal},
  title       = {An Introduction to Delay Differential Equations with Applications to the Life Sciences},
  date        = {2011},
  shortseries = {Texts Appl. Math.},
  series      = {Texts in Applied Mathematics},
  number      = {57},
  publisher   = {Springer},
  location    = {New York},
  doi         = {10.1007/978-1-4419-7646-8},
}

@Manual{Szalai2013,
  author = {Szalai, Róbert},
  title  = {Knut: a continuation and bifurcation software for delay-differential equations (version 8)},
  date   = {2013-12-26},
  url    = {https://rs1909.github.io/knut/},
}

@Book{Trefethen2000,
  author      = {Trefethen, Lloyd Nicholas},
  title       = {Spectral Methods in MATLAB},
  date        = {2000},
  shortseries = {Software Environ. Tools},
  series      = {Software, Environments and Tools},
  number      = {10},
  publisher   = {Society for Industrial and Applied Mathematics},
  location    = {Philadelphia},
  doi         = {10.1137/1.9780898719598},
}

@Article{Trefethen2008,
  author       = {Trefethen, Lloyd Nicholas},
  title        = {Is {G}auss Quadrature Better than {C}lenshaw--{C}urtis?},
  journaltitle = {SIAM Review},
  shortjournal = {SIAM Rev.},
  date         = {2008},
  volume       = {50},
  number       = {1},
  pages        = {67-87},
  doi          = {10.1137/060659831},
}

@Article{VisserMeijerVanPuttenVanGils2012,
  author       = {Visser, Sid and Meijer, Hil G. E. and van Putten, Michel J. A. M. and van Gils, Stephan A.},
  title        = {Analysis of stability and bifurcations of fixed points and periodic solutions of a lumped model of neocortex with two delays},
  journaltitle = {Journal of Mathematical Neuroscience},
  shortjournal = {J. Math. Neurosci.},
  volume       = {2},
  date         = {2012-04-25},
  eid          = {8},
  doi          = {10.1186/2190-8567-2-8},
}

\appendix

%%%%% END FROM ACM TEMPLATE

\section{Pseudospectral discretisation of delay equations}
\label{sec:discretization}

We here summarise the pseudospectral collocation method to approximate a DE (either DDE or RE) with an ODE. For more details, see \cite{BredaDiekmannGyllenbergScarabelVermiglio2016,%
DiekmannScarabelVermiglio2020} for DDEs and \cite{ScarabelDiekmannVermiglio2021} for REs. 

We consider DEs with a finite maximum delay $\tau>0$. 
Given a function $x \colon \mathbb{R} \to \mathbb{R}^d$ for some $d \in \mathbb{N}\setminus\{0\}$, we define \emph{history function}, or \emph{state}, at time $t$, the function
\begin{equation*}
    x_t(\theta)\coloneqq x(t+\theta), \quad \theta \in [-\tau,0].
\end{equation*}

We consider the DDE
\begin{equation}
\label{eq:dde}
x'(t) = F(x_t), \quad t\geq 0,
\end{equation}
with $F \colon X \to \mathbb{R}^d$ and $X \coloneqq C([-\tau,0]; \mathbb{R}^d)$, and the RE
\begin{equation}
\label{eq:re}
x(t) = F(x_t), \quad t>0,
\end{equation}
with $F \colon X \to \mathbb{R}^d$ and $X \coloneqq L^{1}([-\tau,0]; \mathbb{R}^d)$.
We also consider systems that couple DDEs and REs.
We use the same notation in both cases to emphasise the analogy between their treatment.

Using the framework of sun-star calculus \cite{DiekmannGettoGyllenberg2008,DiekmannVanGilsVerduynLunelWalther1995}, each of \eqref{eq:dde} and \eqref{eq:re} can be reformulated as the semilinear abstract differential equation (ADE)
\begin{equation}
\label{eq:ade}
v'(t) = \mathcal{A} v(t) + \mathcal{F}(v(t)), \quad t\geq 0,
\end{equation}
with $v(t) =j x_t$, where:
\begin{itemize}[noitemsep]
    \item $j \colon X \to Y$ is an immersion of the state space $X$ into a suitable space $Y$,
    \item $\mathcal{F}$ is a (nonlinear) perturbation that captures the (nonlinear) action of $F$, and
    \item $\mathcal{A}$ is the infinitesimal generator of the family of solution operators associated to the trivial problem (i.e.\ with $F\equiv 0$), defined as $\mathcal{A}(\psi) = \psi'$ for $\psi$ in a suitable domain $D(\mathcal{A}) \subset Y$.
\end{itemize}
Observe that $\mathcal{F}$ describes the extension of the state into the future, while $\mathcal{A}$ describes the shift of the history.
For DDEs, the `larger' space is $Y\coloneqq\mathbb{R}^d\times L^\infty([-\tau,0]; \mathbb{R}^d)$ and the immersion is $j(\psi)\coloneqq(\psi(0),\psi)$.
For REs, the `larger' space is $Y\coloneqq NBV([-\tau,0]; \mathbb{R}^d)$ and $j$ is the integral operator defined by
\begin{equation}\label{eq:j-re}
(j(\psi))(\theta) \coloneqq -\int_{\theta}^{0} \psi(s) \D s, \quad \theta \in [-\tau, 0],
\end{equation}
so that $v(t)$ is a primitive of the state $x_t$ that takes  value $0$ at $\theta = 0$.

Via pseudospectral collocation, the dynamical system described by the ADE \eqref{eq:ade} in $Y$ is approximated via a finite-dimensional dynamical system by means of suitable restriction and interpolation operators. 

Let $M$ be a positive integer, describing the degree of collocation polynomials, and choose a family of nodes $\{\theta_{k}\}_{k \in \{0, \dots, M\}}$ in the interval $[-\tau, 0]$ such that 
\begin{equation*}
    0 = \theta_{0} > \theta_{1} > \dots > \theta_{M} \geq -\tau. 
\end{equation*}
A continuous function $\psi$ with domain $[-\tau, 0]$ is projected to a finite-dimensional space $\mathbb{R}^{d(M+1)}$ via restriction on the nodes:
\begin{equation*}
    \mathcal{R}_M \psi = (\psi(\theta_0);\dots; \psi(\theta_M)) \in \mathbb{R}^{d(M+1)},
\end{equation*}
where the semicolon concatenates vectors vertically. 
Vice versa, given the vectors $\Psi_0 \in \mathbb{R}^d$ and $\Psi = (\Psi_1;\dots;\Psi_M) \in \mathbb{R}^{dM}$ with $\Psi_k \in \mathbb{R}^d$, $k\in\{1,\dots,M\}$, we define an interpolation operator $\mathcal{P}_M$ that maps $(\Psi_0;\Psi)$ to the interpolating $d$-valued polynomial of degree at most $M$, namely
\begin{equation*}
    \mathcal{P}_M((\Psi_0;\Psi)) = \sum_{k=0}^M \Psi_k \ell_k,
\end{equation*}
where 
\begin{equation*}
    \ell_k(\theta) =  \prod_{\substack{i=0\\ i\neq k}}^M \frac{\theta-\theta_i}{\theta_k-\theta_i}, \quad k \in \{0,\dots,M\},
\end{equation*}
are the Lagrange polynomials corresponding to the collocation points $\{\theta_k\}_{k \in \{0, \dots, M\}}$. 

To write the approximating system in a compact way, we also define the differentiation matrix $\widehat{D}_{M}\in\mathbb{R}^{(M+1)\times(M+1)}$ associated with the nodes, such that 
\begin{equation*}
    (\widehat{D}_M)_{ij} = \ell_j'(\theta_i), \quad i,j \in\{0,\dots,M\}.
\end{equation*}
We denote by $D_{M}$ the matrix obtained from $\widehat{D}_{M}$ removing both the first row and the first column.
We also denote the matrices obtained from $\widehat{D}_{M}$ removing either the first row or the first column as, respectively, $D_{M,c}$ (`$D_{M}$ plus the first column') and $D_{M,r}$ (`$D_{M}$ plus the first row').

In the case of the DDE \eqref{eq:dde}, the equivalent ADE \eqref{eq:ade} is discretised as the $d(M+1)$-dimensional ODE
\begin{equation}
\label{eq:dde-discr}
\left\{
\begin{aligned}
&x'_{M} = F(\mathcal{P}_{M}((x_M; V_M))), \\
&V'_{M} = (D_{M,c} \otimes I_{d}) (x_{M}; V_{M}),
\end{aligned}
\right.
\end{equation}
for $x_M \in \mathbb{R}^d$ and $V_M = (V_{M,1};\dots;V_{M,M}) \in \mathbb{R}^{dM}$, where $I_{d} \in \mathbb{R}^{d\times d}$ is the identity matrix and $\otimes$ is the Kronecker tensor product.
For all $t\geq 0$, the state $v(t) = jx_t = (x(t),x_t)$ is approximated by the pair $(x_M(t),\mathcal{P}_M((x_M(t);V_M(t)))$, so that $x_M(t) \in \mathbb{R}^d$ can be taken as an approximation of the solution $x(t)$ and each $V_{M,k}$ as an approximation of the translation $x(t+\theta_k)$, for $k\in\{1,\dots,M\}$. 
If $\psi \in C([-\tau,0],\mathbb{R}^d)$ is an initial condition (or equilibrium) for \eqref{eq:dde}, the corresponding initial condition (or equilibrium) for \eqref{eq:dde-discr} is obtained by taking the restriction $\mathcal{R}_M \psi \in \mathbb{R}^{d(M+1)}$.

For the RE \eqref{eq:re}, the equivalent ADE \eqref{eq:ade} is discretised as the $dM$-dimensional ODE
\begin{equation}
\label{eq:re-discr}
V'_{M} = (D_{M,c} \otimes I_{d}) (0_d; V_{M}) - F_M(\mathcal{P}_{M} ((D_{M,r} \otimes I_d) V_{M})),
\end{equation}
for $V_M \in \mathbb{R}^{dM}$, 
where
$F_M(\cdot)=(F(\cdot);\dots;F(\cdot))$, with $F(\cdot)$ repeated $M$ times,
and $0_d \in\mathbb{R}^{d}$ is the zero vector.
The state $v(t) = jx_t$, with $j$ defined by \eqref{eq:j-re}, is approximated by the polynomial
$\mathcal{P}_M((D_{M,r} \otimes I_d) V_{M}))$.
In fact, since $v(t)(0)$ is always $0$, there is no component of the ODE describing the integrated state at $\theta_0=0$.
For each time $t$, the solution $x(t)$
is approximated by
\begin{equation}
\label{eq:re-extra-step}
x_{M}(t)\coloneqq F(\mathcal{P}_{M} ((D_{M,r} \otimes I_d) V_{M})) \in \mathbb{R}^{d}.
\end{equation}
If $\psi \in L^1([-\tau,0],\mathbb{R}^d)$ is an initial condition (or equilibrium) for \eqref{eq:re}, the corresponding initial condition (or equilibrium) for \eqref{eq:re-discr} is given by the vector $\mathcal{R}_M j\psi$ when the first block of $d$ components (relevant to the node $\theta_0=0$) is removed; the resulting vector $\Psi\in\mathbb{R}^{dM}$ can be approximated by
\begin{equation}
\label{eq:re-initial-value}
    (D_M^{-1}\otimes I_d) (\psi(\theta_1);\dots;\psi(\theta_M)) \in\mathbb{R}^{dM}.
\end{equation}
Note that, when $\psi \equiv \overline{\psi}$ is an equilibrium of \eqref{eq:re}, then $\mathcal{R}_M j\psi = (\overline{\psi}\theta_1; \dots; \overline{\psi}\theta_M)$.

Finally, if distributed delays are present, a quadrature formula needs to be applied in order to fully discretise the equation.

In the implementation, as collocation nodes we use the family of Chebyshev points of the second kind, rescaled to the maximal delay interval $[-\tau,0]$:
\begin{equation*}
    \theta_k = \frac{\tau}{2}\left(\cos\left(\frac{k\pi}{M}\right)-1\right), \quad k \in\{0,\dots,M\}.
\end{equation*}
Approximation methods based on them typically exhibit spectral accuracy in the order $M$ of the approximation (see, e.g., \cite{Trefethen2000,MastroianniMilovanovic2008} for more details).
Moreover, these nodes include the endpoints of the interval, which allows to avoid the interpolation for the values of the coordinates at the maximum delay for DDEs in many practical cases.
For these nodes, the matrix $D_M^{-1}$ in \eqref{eq:re-initial-value} approximates the integral operator $j$ as desired \cite{DiekmannScarabelVermiglio2020}.

For the quadrature of the integrals representing distributed delays we use the Clenshaw--Curtis formula \cite{ClenshawCurtis1960,Trefethen2008}, which is again based on the chosen Chebyshev points, rescaled in the integration interval (see \cite[chapter 12]{Trefethen2000} for the computation of the quadrature weights).

Interpolation is done with the barycentric Lagrange interpolation formula \cite{BerrutTrefethen2004}.
We also note that we avoided the use of \texttt{kron} where possible in view of efficiency.

%%%%% FROM ACM TEMPLATE

\end{document}